%% file: main.tex
\documentclass{aero}

\input{sections/preamble.tex}

\ADsetup{
title    = {Non-linear stochastic trajectory optimisation},
author   = {Thomas Caleb$^{1}$\cor{}, Roberto Armellin$^{2}$, and St\'ephanie Lizy-Destrez$^{1}$},
address  = {
	1. DCAS/SaCLaB, ISAE-SUPAERO, Toulouse 31055, France \sep
	2. Te P\=unaha \=Atea -- Space Institute, University of Auckland, Auckland 1010, New Zealand
},
email    = {Thomas Caleb, thomas.caleb@dycsyt.com},
abstract = {Designing robust space trajectories in non-linear dynamical environments, such as the Earth–Moon \gls*{CR3BP}, poses significant challenges due to sensitivity to initial conditions and non-Gaussian uncertainty propagation. This work introduces a novel solver for discrete-time chance-constrained trajectory optimisation under uncertainty, referred to as \gls*{SODA}. \gls*{SODA} combines \gls*{DA} with adaptive Gaussian mixture decomposition to efficiently propagate non-Gaussian uncertainties, and enforces Gaussian multidimensional chance constraints. 
This work further introduce a risk allocation strategy across mixture components that enables tight and adaptive distribution of safety margins.
The framework is validated on four trajectory design problems of increasing dynamical complexity, from heliocentric transfers to challenging Earth–Moon \gls*{CR3BP} scenarios. A linear variant, the \gls*{L-SODA} solver, recovers deterministic performance with minimal overhead under small uncertainties, while the non-linear \gls*{SODA} solver yields improved robustness and tighter constraint satisfaction in strongly non-linear regimes. 
Results highlight \gls*{SODA}’s ability to generate accurate, robust, and computationally tractable solutions, supporting its potential for future use in uncertainty-aware space mission design.},
keywords = {Stochastic trajectory optimisation \sep Non-linear dynamics \sep High-order Taylor expansions \sep Chance constraints \sep Gaussian mixture models}
}

\begin{document}

\maketitle

% Sections
% \input{sections/0_abstract.tex}
\input{sections/1_introduction.tex}
\input{sections/2_background.tex}
\input{sections/3_l-soda}
\input{sections/4_soda}
\input{sections/5_application}
\input{sections/6_conclusion}
\input{sections/7_funding_acknoledgments}
% \appendix
% \input{sections/8_appendix}

% Bib
\section*{References}
\bibliographystyle{unsrtnat}
\bibliography{references}

\end{document}

%% file: sections/preamble.tex
\usepackage[utf8]{inputenc}
\usepackage{textcomp}

\usepackage[version=4]{mhchem}
\usepackage{siunitx}
\usepackage{lipsum}

\usepackage[nonumberlist,acronym]{glossaries} % enable acronyms
\usepackage{graphicx,setspace,parskip,bm,amsmath,mathrsfs} % enable figures, symbols, etc.

\usepackage{amsthm}
\usepackage{longtable,tabularx,multicol,multirow,threeparttablex} % enable tables
\usepackage{float} % enable forcing figures/tables placement
\usepackage{makecell} % to customize cells in table
\usepackage{booktabs} % enable rule, toprule, midrule, etc.
\usepackage{upgreek}
\usepackage{footnpag,marginnote} % enable footnotes and side notes
\usepackage{booktabs} % enable rule, toprule, midrule, etc.
\usepackage[font=small]{caption} % enable caption
\usepackage{subcaption} % enable subfigure and subcaption
\usepackage{comment}
\usepackage{microtype} % slightly tweak font spacing for aesthetics
\usepackage{makecell} % to customize cells in table
\usepackage{setspace} % enable setting space size
\usepackage{pifont}
\usepackage{xcolor} % required for specifying colors by name
\usepackage{colortbl} % add colors to tables
\usepackage{mathtools}
\usepackage{centernot}
\setlist{labelindent=\parindent,leftmargin=*,nosep}
\setlist[enumerate]{label={\arabic*)}}
\usepackage{pgfplots}
\pgfplotsset{compat=newest}
\pgfplotsset{plot coordinates/math parser=false}
\newlength\figureheight
\newlength\figurewidth
\usepackage{environ,tikz} % enable tikz figures
\usetikzlibrary{mindmap, trees, calc}
\usetikzlibrary{positioning,fit}
\usetikzlibrary{through}
\usetikzlibrary{shapes.geometric, arrows}
\usetikzlibrary{automata,positioning}

\tikzset{block/.style={draw,thick,text width=2cm,minimum height=1cm,align=center},
         line/.style={-latex}
}
\usepackage{etoolbox}
\usepackage{algorithm}
\usepackage{algpseudocode}

\newtheorem{theorem}{Theorem}

\DeclareSIUnit{\AU}{AU}
\DeclareSIUnit{\day}{d}
\DeclareSIUnit{\deg}{deg}
\DeclareSIUnit{\year}{yr}

\newcommand{\Eq}[1]{Eq.~\eqref{#1}}

\newcommand{\Fig}[1]{Fig.~\ref{#1}}

\newcommand{\Tab}[1]{Table~\ref{#1}}

\newcommand{\Thm}[1]{Theorem~\ref{#1}}

\newcommand{\Alg}[1]{Alg.~\ref{#1}}
\newcommand{\Sec}[1]{Section~\ref{#1}}

\newcommand{\ie}{i.\,e.,~}

\newcommand{\eg}{e.\,g.,~}

\newcommand{\lastdate}{May 1, 2026}

\newacronym{ABM}{ABM}{Adams-Bashforth-Moulton}
\newacronym{ADS}{ADS}{automatic domain splitting}
\newacronym{AUL}{AUL}{augmented Lagrangian}
\newacronym{CDF}{CDF}{cumulative distribution function}
\newacronym{CR3BP}{CR3BP}{circular restricted three-body problem}
\newacronym{EMS}{EMS}{Earth—Moon system}
\newacronym{DRO}{DRO}{distant retrograde orbit}
\newacronym{DROs}{DROs}{distant retrograde orbits}
\newacronym{DA}{DA}{differential algebra}
\newacronym{DACE}{DACE}{differential algebra core engine}
\newacronym{DADDy}{DADDy}{differential algebra-based differential dynamic programming}
\newacronym{DDP}{DDP}{differential dynamic programming}
\newacronym{DNC}{DNC}{did not converge}
\newacronym{DOPRI8}{DOPRI8}{Dormand-Prince 8th-order embedded Runge-Kutta method}
\newacronym{DOP853}{DOP853}{Dormand-Prince of order 8(5,3) embedded Runge-Kutta method}
\newacronym{EOM}{EoM}{equations of motion}
\newacronym{ESA}{ESA}{European Space Agency}
\newacronym{GEO}{GEO}{geostationary orbit}
\newacronym{GMM}{GMM}{Gaussian mixture model}
\newacronym{GTO}{GTO}{geostationary transfer orbit}
\newacronym{IC}{IC}{initial conditions}
\newacronym{iLQR}{iLQR}{iterative linear-quadratic regulator}
\newacronym{JAXA}{JAXA}{Japan aerospace exploration agency}
\newacronym{JPL}{JPL}{Jet Propulsion Laboratory}
\newacronym{MC}{MC}{Monte-Carlo}
\newacronym{MB}{MB}{megabytes}
\newacronym{L-SODA}{L-SODA}{linear stochastic optimisation with differential algebra}
\newacronym{LAM}{LAM}{likelihood agreement measure}
\newacronym{LEO}{LEO}{low Earth orbit}
\newacronym{LOADS}{LOADS}{low-order automatic domain splitting}
\newacronym{LQG}{LQG}{linear quadratic Gaussian}
\newacronym{MEO}{MEO}{medium Earth orbit}
\newacronym{NASA}{NASA}{National Aeronautics and Space Administration}
\newacronym{NEA}{NEA}{near-Earth asteroid}
\newacronym{NEO}{NEO}{near-Earth object}
\newacronym{NLI}{NLI}{non-linearity index}
\newacronym{NRHO}{NRHO}{near-rectilinear halo orbit}
\newacronym{NRHOs}{NRHOs}{near-rectilinear halo orbits}
\newacronym{NSG}{NSG}{non-spherical gravity}
\newacronym{ODE}{ODE}{ordinary differential equation}
\newacronym{PECE}{PECE}{predictor-corrector}
\newacronym{PCK}{PCK}{planetary constants kernel}
\newacronym{PDE}{PDE}{partial differential equation}
\newacronym{PDF}{PDF}{probability density function}
\newacronym{PDS}{PDS}{Planetary Data System}
\newacronym{PN}{PN}{projected Newton}
\newacronym{RAAN}{RAAN}{right ascension of the ascending node}
\newacronym{RHS}{RHS}{right-hand side}
\newacronym{RK}{RK}{Runge-Kutta}
\newacronym{RK78}{RK78}{order 7(8) Runge–Kutta–Fehlberg embedded method}
\newacronym{RSS}{RSS}{root sum square}
\newacronym{RTN}{RTN}{radial-tangential-normal}
\newacronym{RPF}{RPF}{roto-pulsating frame}
\newacronym{SDE}{SDE}{stochastic differential equation}
\newacronym{SODA}{SODA}{stochastic optimisation with differential algebra}
\newacronym{STT}{STT}{state transition tensor}
\newacronym{STTs}{STTs}{state transition tensors}
\newacronym{STM}{STM}{state-transition matrix}
\newacronym{TDA}{TDA}{Taylor Differential Algebra}
\newacronym{ToF}{ToF}{time-of-flight}
\newacronym{TPBVP}{TPBVP}{two-point boundary value problem}
\newacronym{UOM}{uom}{unit of measurement}
\newacronym{UT}{UT}{unscented transform}
\newacronym{VSVO}{VSVO}{variable-step variable-order}
\newacronym{VV}{V\&V}{verification and validation}

%% file: sections/1_introduction.tex
\section{Introduction}
\label{sec:introduction}
\glsresetall
Interplanetary trajectory design has grown in complexity, frequently requiring the use of multi-body dynamics and the exploitation of natural dynamical structures. Recent missions such as \gls*{ESA}'s \textit{JUICE} \citep{BoutonnetEtAl_2024_DtJT} have made extensive use of planetary swing-bys and ballistic capture \citep{TopputoBelbruno_2015_ETwBC}, while the \gls*{CR3BP} continues to play a central role in cis-lunar navigation \citep{SmithEtAl_2020_TAPaOoNAtRHttM}. These non-linear dynamical environments are highly sensitive to initial conditions, due to their inherent chaotic nature \citep{Poincare_1892_LMNDLMC}, making trajectory robustness particularly difficult to assess and enforce. Small deviations in state or control, caused by navigation errors or imperfect manoeuvres, can lead to significant divergences from the nominal path. As a result, traditional deterministic optimisation methods often fall short in guaranteeing constraint satisfaction for high-stakes missions.

Traditionally, robustness is addressed \textit{a posteriori} using Monte-Carlo simulations and conservative design margins. This approach mitigates risk but produces over-cautious solutions and does not guarantee \textit{a priori} constraint satisfaction. Stochastic optimal control offers an alternative by managing uncertainty during optimisation \citep{Zakai_1969_OtOFoDP}.
Theoretically, stochastic optimal control enables direct propagation of probability distributions. However, solving the associated partial differential equations (e.g., Fokker–Planck or Hamilton–Jacobi–Bellman) \citep{YongZhou_1999_DPaHE} becomes intractable in high-dimensional or constrained settings.

Recent stochastic optimisation solvers for mission analysis typically follow a modular three-block architecture. First, uncertainty propagation techniques allow to compute the evolution of the state and control distributions. Methods include linear covariance propagation used in \citet{BenedikterEtAl_2022_CAtCCwAtSLTTO}, unscented transforms \citep{JulierUhlmann_1997_NEotKFtNS} in \citet{OzakiEtAl_2020_TSOCfNCTOP}, \gls*{GMM} \citep{DeMarsEtAl_2013_EBAfUPoNDS} in \citet{BooneMcMahon_2022_NGCCTCUGMaRA}, and \gls*{MC} \citep{RobertCasella_2004_MCSM} sampling in \citet{BlackmoreEtAl_2010_APPCAoCCSPC}. A recent advance, the \gls*{LOADS} method from \citet{LosaccoEtAl_2024_LOADSAfNUM}, offers scalable solutions for describing non-Gaussian distributions using Taylor polynomials.

Second, a transcription method converts probabilistic constraints and objective functions into tractable deterministic sufficient conditions. This step, often overlooked, can introduce significant conservatism. While efficient methods exist for one-dimensional cases or specific constraint structures \citep{BlackmoreEtAl_2011_CCOPPwO, RidderhofEtAl_2020_CCCCfLTMFTO, OguriLantoine_2022_SSCPfRLTTDuU}, general-purpose approaches for Gaussian chance constraints remain scarce. Moreover, existing transcription techniques typically require eigenvalue computations, resulting in cubic complexity. Reducing this computational burden enables the resolution of more ambitious problems and broaden the applicability of such methods. Our recent first-order transcription method \citep{CalebEtAl_2025_CCTaFREfSTO} provides a practical alternative: it handles multidimensional Gaussian constraints with low conservatism and linear complexity.

Third, the resulting transcribed problem is formulated, e.g., as a non-linear programming problem \citep{GrecoEtAl_2022_RSTDUBOC,MarmoEtAl_2023_AHMSAfCCoIMwNE}, or a convex problem \citep{RidderhofEtAl_2020_CCCCfLTMFTO,BenedikterEtAl_2022_CAtCCwAtSLTTO}. Despite this architecture, several challenges persist. Many methods rely on strong assumptions such as gaussianity or linear dynamics \citep{OzakiEtAl_2018_SDDPwUTfLTTD,BenedikterEtAl_2022_CAtCCwAtSLTTO}, or exhibit high computational cost \citep{GrecoEtAl_2022_RSTDUBOC}.
Moreover, transcription methods often remain overly conservative \citep{RidderhofEtAl_2020_CCCCfLTMFTO}. To mitigate these issues, some recent works develop \textit{ad hoc} solvers, such as stochastic variants of \gls*{DDP} \citep{OzakiEtAl_2018_SDDPwUTfLTTD} or methods based on primer vector theory \citep{OguriMcMahon_2022_SPVfRLTTDuU}. However, the field still lacks general-purpose, unified frameworks capable of coupling uncertainty propagation, constraint transcription, and optimisation within a coherent and computationally efficient pipeline.

In this work, we present the \gls*{SODA} solver, a unified stochastic optimisation framework that directly addresses these limitations. Most existing methods for robust trajectory optimisation either assume Gaussian statistics or rely on conservative approximations to ensure tractability. Such approximations, particularly at the transcription stage, often lead to overly conservative designs and limit the effectiveness of robustness.
Our first contribution is a second-order discrete-time \gls*{DDP} framework, originally introduced in \citet{CalebEtAl_2025_APBCSfFOLTTO} and extended here to the stochastic setting. Within this framework, uncertainty propagation, dynamics evaluation, constraint transcription, and automatic differentiation are performed within a single efficient computational pipeline based on high-order Taylor expansions.

We then leverage this solver as a building block to address non-Gaussian uncertainty using an adaptive \gls*{GMM}. To this end, we develop a scalable uncertainty propagation method that combines \gls*{DA} with adaptive decomposition, extending the \gls*{LOADS} algorithm \citep{LosaccoEtAl_2024_LOADSAfNUM} to capture non-Gaussian effects with high fidelity. This integrated use of \gls*{DA} reduces numerical overhead and improves scalability to large-scale problems.

The effectiveness of the proposed approach is demonstrated through four representative trajectory design test cases, covering both heliocentric and Earth--Moon \gls*{CR3BP} environments. The results show that \gls*{SODA} produces robust trajectories that satisfy chance constraints even in highly non-linear regimes, while achieving performance close to the deterministic baseline within practical computational times.

The remainder of this article is organised as follows. In \Sec{sec:background}, we provide the necessary background on constrained \gls*{DDP}, chance-constrained optimal control problems, transcription methods and risk estimation, \gls*{GMM}, and the \gls*{DA} framework.
In \Sec{sec:LSODA}, we introduce our \gls*{L-SODA} solver, which serves as a building block for the \gls*{SODA} solver, presented in \Sec{sec:soda}. \Sec{sec:application} demonstrates the application and validation of \gls*{L-SODA} and \gls*{SODA} through a series of four trajectory optimisation problems in the two-body problem and the Earth--Moon \gls*{CR3BP}. Finally, we summarise our conclusions in \Sec{sec:conclusion}.

%% file: sections/2_background.tex
\section{Background}
\label{sec:background}
\glsresetall
\subsection{Constrained differential dynamic programming}
\subsubsection{Differential dynamic programming}
\Gls*{DDP} addresses finite-horizon optimal control problems governed by discrete-time dynamics of the form $\bm{x}_{k+1} = \bm{f}(\bm{x}_k, \bm{u}_k)$, where $k \in \mathbb{N}_{N-1}=\{0,1,\cdots,N-1\}$, $\bm{x}_k \in \mathbb{R}^{N_x}$ is the state, and $\bm{u}_k \in \mathbb{R}^{N_u}$ is the control input \citep{Mayne_1966_ASOGMfDOToNLDTS,LantoineRussell_2012_AHDDPAfCOCPP1T}. Given the initial state $\bm{x}_0$ and a target state $\bm{x}_t$, the goal is to minimise a cost function: $J(\bm{U}) = \sum_{k=0}^{N-1} l(\bm{x}_k, \bm{u}_k) + \phi(\bm{x}_N, \bm{x}_t)$, where $\bm{U} = (\bm{u}_0, \cdots, \bm{u}_{N-1})$, $l$ is the stage cost, and $\phi$ is the terminal cost.
Using Bellman's principle of optimality, the value function $V_k(\bm{x}_k)$ defines the optimal cost-to-go from step $k$, leading to the recursive definition of the $Q$-function: $Q_k(\bm{x}_k, \bm{u}_k) = l(\bm{x}_k, \bm{u}_k) + V_{k+1}(\bm{f}(\bm{x}_k, \bm{u}_k))$. The partial derivatives of $Q_k$ are used to compute the optimal control policy.

Second-order terms in the $Q$ function derivatives account for the curvature of both the cost function and the system dynamics. In particular, \gls*{DDP} includes second-order derivatives of $\bm{f}$, in contrast to \gls*{iLQR}, which neglects them to decrease run time.
These terms enhance the local approximation of the value function and enable quadratic convergence \citep{NgangaWensing_2021_ASODDPfRBS}, although they often require regularisation for numerical robustness \citep{LantoineRussell_2012_AHDDPAfCOCPP1T}. 
The trajectory $\bm{X} = (\bm{x}_0, \cdots, \bm{x}_N)$ and the control $\bm{U}$ are updated by applying the control corrections $\bm{a}_k$ (feed-forward) and $\bm{K}_k$ (feedback), propagated forward through the system dynamics. The updated states and controls, denoted $\bm{x}^*_k$ and $\bm{u}^*_k$, are computed iteratively. The initial state remains unchanged, i.e., $\bm{x}^*_0 = \bm{x}_0$. They are retrieved iteratively using the corrections $\delta\bm{x}^*_k$ and $\delta\bm{u}^*_k$: $\delta\bm{x}^*_k = \bm{x}_k^* - \bm{x}_k$, $\delta\bm{u}^*_k = \bm{a}_k + \bm{K}_k \delta\bm{x}^*_k$, $\bm{u}_k^* = \bm{u}_k + \delta\bm{u}^*_k$, and $\bm{x}_{k+1}^* = \bm{f}(\bm{x}_k^*, \bm{u}_k^*)$.
At convergence, the feed-forward terms satisfy $\|\bm{a}_k\| \approx 0$, indicating that the trajectory has reached a local minimum. The feedback gains $\bm{K}_k$ can then be used to reject disturbances or recover from deviations due to modelling or navigation errors.

\subsubsection{Constraints handling}
The original \gls*{DDP} formulation does not handle constraints \citep{Mayne_1966_ASOGMfDOToNLDTS}. Such constraints are essential in practical applications, for instance, to enforce actuator limits or ensure safety. Following the work of \citet{HowellEtAl_2019_AaFSfCTO}, this limitation is addressed through an augmented Lagrangian framework. This framework accommodates path constraints of the form $\bm{g}(\bm{x}_k, \bm{u}_k) \preceq \bm{0}$ of size $N_g$ and terminal constraints of the form $\bm{g}_{t}(\bm{x}_N, \bm{x}_t) \preceq \bm{0}$ of size $N_{gt}$. 
We define $\bm{y}_1\preceq \bm{y}_2$ such that for all $i$ the $i$\textendash{th} component of $\bm{y}_1$ is less than or equal to the $i$\textendash{th} component of $\bm{y}_2$.
This approach allows the solver to handle a broad class of constrained trajectory optimisation problems while retaining the efficiency and structure of the original \gls*{DDP} method.

\subsubsection{Differential algebra-based constrained \gls*{DDP} solver}

\citet{CalebEtAl_2025_APBCSfFOLTTO} accelerate \gls*{DDP} by leveraging the \gls*{DA} framework \citep{Berz_1999_MMMiPBP}, which systematically manipulates high-order Taylor expansions of multivariate functions. This enables both automatic differentiation and local polynomial approximation of non-linear dynamics.
In \gls*{DA}, a function $\bm{h}$ is expanded around a nominal point $\bm{z}$ as a polynomial $\mathcal{P}_{\bm{h}}(\delta \bm{z})$.
In this paper, the constant part of a polynomial $\mathcal{P}_{\bm{h}}$ is denoted $\overline{\mathcal{P}_{\bm{h}}}$, a small perturbation of any given quantity $z$ is denoted $\delta z$, and the derivative of a function $\bm{h}$ with respect to $\bm{z}$ is denoted $\bm{h}_{\bm{z}}$.
Derivatives are extracted directly from the expansion, while the approximation is valid within a convergence radius $R_\varepsilon$, satisfying \citet{WittigEtAl_2015_PoLUSiODbADS}: $\|\delta\bm{z}\| \leq R_{\varepsilon} \implies \left\|\mathcal{P}_{\bm{h}}(\delta\bm{z}) - \bm{h}(\bm{z}+\delta\bm{z})\right\| \leq \varepsilon$.
This approach is especially effective when function evaluations are costly or repeated, as in trajectory or uncertainty propagation \citep{ArmellinEtAl_2010_ACECUDAtCoA,WittigEtAl_2015_PoLUSiODbADS,CalebEtAl_2023_DAMAtCAGaBDAtPFotES}. Within \gls*{DDP}, it accelerates the forward and backward passes by replacing standard evaluations and derivatives with polynomial-based ones, provided the corrections remain within $R_{\varepsilon_{\textrm{DA}}}$.
The resulting \gls*{DADDy} solver\footnote{Available at: \url{https://github.com/ThomasClb/DADDy.git} [last accessed \lastdate].}, improves runtime without compromising precision or convergence. It also supports constraint enforcement via the augmented Lagrangian strategy described previously. The complete method is detailed in \citet{CalebEtAl_2025_APBCSfFOLTTO}.

\subsubsection{\Gls*{DA}-based Newton solver}
While the \gls*{AUL} formulation enforces constraints up to a tolerance $\varepsilon_{\textrm{AUL}}$, \citet{HowellEtAl_2019_AaFSfCTO} propose using a Newton method to achieve higher accuracy, with a tighter tolerance $\varepsilon_{\textrm{N}} < \varepsilon_{\textrm{AUL}}$. 
It consists of concatenating all the active constraints in a single vector $\bm{\Gamma}$ of size $N_\Gamma \leq N\cdot (N_{g} + N_x) + N_{gt} $, including continuity equality constraints: $\bm{g}_{\textrm{c},k}=\bm{g}_{\textrm{c}}\left(\bm{x}_k, \bm{u}_k, \bm{x}_{k+1}\right)=\bm{x}_{k+1} -  \bm{f}\left(\bm{x}_k, \bm{u}_k\right)=\bm{0}$. The controls and the state vectors are concatenated in a vector $\bm{Y}$ of size $N_Y=N\cdot(N_x+N_u)$, defined as: $\bm{Y}=\left[\bm{u}_0^{{\textrm{T}}},\bm{x}_1^{{\textrm{T}}},\bm{u}_0^{{\textrm{T}}}\cdots,\bm{x}_{N-1}^{{\textrm{T}}},\bm{u}_{N-1}^{{\textrm{T}}}, \bm{x}_{N}^{{\textrm{T}}}\right]^{{\textrm{T}}}$.
The active constraints are linearised as: $\bm{\Gamma}\left(\bm{Y}+\delta\bm{Y}\right)\approx \bm{\Delta}\delta\bm{Y}+\Bar{\bm{\Gamma}}$, where $\bm{\Delta}$ and $\Bar{\bm{\Gamma}}$ denote, respectively, the gradient and the value of the constraints evaluated at $\bm{Y}$, and $\delta\bm{Y}$ is a small variation around $\bm{Y}$.
The correction $\delta\bm{Y}^* = -\bm{\Delta}^{+} \Bar{\bm{\Gamma}}$, where $\bm{\Delta}^{+}=\bm{\Delta}^{\textrm{T}}\left(\bm{\Delta}\bm{\Delta}^{\textrm{T}}\right)^{-1}$, yields a refined estimate $\bm{Y}^*$.
To further accelerate the Newton solver, \citet{CalebEtAl_2025_APBCSfFOLTTO} apply the \gls*{DA} framework for automatic differentiation and dynamics evaluation, similar to the \gls*{DA}-based \gls*{DDP} solver. In addition, due to the symmetric block tri-diagonal structure of $\bm{\Delta}\bm{\Delta}^{\textrm{T}}$, the computation of $\delta\bm{Y}^*$ can be accelerated using block Cholesky factorisation \citep{Cholesky_1910_SLRNDSDL,CaoEtAl_2002_PCFoaBTM}, bringing the complexity from cubic to quadratic. 

\subsection{Automatic Gaussian mixture model decomposition of a Gaussian distribution}
\label{sec:GMM}
A key challenge in stochastic control lies in handling complex, non-Gaussian distributions. To address this, an automatic \gls*{GMM} decomposition is employed to reformulate the original problem as a collection of simpler Gaussian sub-problems.

\subsubsection{Gaussian mixture model decomposition of a Gaussian distribution}
A \gls*{GMM} approximates the \gls*{PDF} $p$ of a multivariate random variable $\bm{y}$ using a weighted sum of multivariate Gaussian distributions $p_{\textrm{G}}$: $p\left(\bm{y}\right)\approx p_{\textrm{GMM}}\left(\bm{y};\bm{\alpha},\Bar{\bm{\mu}},\bm{S}\right) = \sum_{i=1}^{M}{\alpha_i\, p_{\textrm{G}}\left(\bm{y};\Bar{\bm{\mu}}_i, \bm{S}_i\right)}$, where $ \bm{\alpha} = \left(\alpha_1, \alpha_2, \cdots, \alpha_M\right) $, $ \Bar{\bm{\mu}} = \left(\Bar{\bm{\mu}}_1, \Bar{\bm{\mu}}_2, \cdots,\Bar{\bm{\mu}}_M\right) $, $ \bm{S} = \left(\bm{S}_1, \bm{S}_2, \cdots, \bm{S}_M\right) $, and $ \sum_{i=1}^{M} \alpha_i = 1 $.
Following the methodology of \citet{DeMarsEtAl_2013_EBAfUPoNDS}, $ p = p_{\textrm{G}}\left(\cdot; 0, 1\right) $, a univariate standard Gaussian distribution, is approximated with a sum of $ M = \num{3} $ components. The variances of the Gaussian terms are assumed equal to a fixed value $ \sigma^2 $, and the first component is centred at zero to simplify the parameter identification. Due to the symmetry of $ p_{\textrm{G}} $, the \gls*{GMM} reduces to:
\begin{equation}
    \label{eq:univariate_GMM}
    p_{\textrm{GMM}}\left(y;\alpha,\Bar{\mu},\sigma^2\right) = \alpha\, p_{\textrm{G}}\left(y;0, \sigma^2\right) + \frac{1-\alpha}{2}\left[ p_{\textrm{G}}\left(y;\Bar{\mu}, \sigma^2\right) + p_{\textrm{G}}\left(y;-\Bar{\mu}, \sigma^2\right)\right].
\end{equation}
The parameters $ (\alpha,\Bar{\mu},\sigma^2) $ are obtained by solving an optimisation problem. The numerical solution of that problem used throughout this article is from \citet{DeMarsEtAl_2013_EBAfUPoNDS}.
\begin{comment}
\begin{table}
    \centering
    \caption{\gls*{GMM} parameters for $ M=3 $ from \citet{DeMarsEtAl_2013_EBAfUPoNDS}.}
    \footnotesize
    \begin{tabular}{c c c} 
        \toprule
        \textbf{$ \alpha $} & \textbf{$ \Bar{\mu} $} & \textbf{$ \sigma $} \\ [0.5ex] 
        \midrule
        \num{0.5495506294920584} & \num{1.0575150485760967} & \num{0.6715664864669252} \\
        \bottomrule
    \end{tabular}
    \small
    \label{tab:GMM_DeMars_Fossa_K3}
\end{table}
\end{comment}
This univariate approximation serves as a foundation for extending the \gls*{GMM} to the multivariate case of a $d$-dimensional Gaussian distribution $\bm{y} \sim \mathcal{N}\left(\Bar{\bm{y}}, \bm{\Sigma}_{\bm{y}}\right)$, with $d \in \mathbb{N}^*$. The covariance matrix $\bm{\Sigma}_{\bm{y}}$ is symmetric and positive definite, and thus can be orthogonally diagonalised as $\bm{\Sigma}_{\bm{y}} = \bm{O}\bm{\Lambda}\bm{O}^{\textrm{T}}$, where $\bm{O}$ is a rotation matrix and $\bm{\Lambda} = \textrm{diag}\left(\lambda_1,\lambda_2,\cdots,\lambda_d\right)$ is diagonal with $\lambda_j > 0$. 

Let $j \in \mathbb{N}_d^*$ be the decomposition direction. \Eq{eq:univariate_GMM} is then extended to the multivariate case as:
\begin{equation}
    \label{eq:mulivariate_GMM}
    p_{\textrm{GMM}}\left(\bm{y};\alpha,\Bar{\bm{\mu}},\bm{S}\right) = \alpha\, p_{\textrm{G}}\left(\bm{y};\Bar{\bm{y}}, \bm{S}\right) + 
    \frac{1-\alpha}{2} \left[p_{\textrm{G}}\left(\bm{y}; \Bar{\bm{y}} -\Bar{\bm{\mu}}, \bm{S}\right)
    + p_{\textrm{G}}\left(\bm{y};\Bar{\bm{y}}+\Bar{\bm{\mu}}, \bm{S}\right)\right].
\end{equation}
The multivariate parameters $ \alpha, \Bar{\bm{\mu}}, \bm{S} $ are derived from the univariate values of \citet{DeMarsEtAl_2013_EBAfUPoNDS} as $ \Bar{\bm{\mu}} = \Bar{\mu} \sqrt{\lambda_j}\, \bm{O}_j$ and $\bm{S} = \bm{O} \, \textrm{diag}\left(\lambda_1,\cdots,\sigma^2\lambda_j,\cdots,\lambda_d\right)\, \bm{O}^{\textrm{T}}$, where $\bm{O}_j$ denotes the $j$-th column of $\bm{O}$.
Note that, since $\sigma < 1$, the iso-\gls*{PDF} ellipsoid described by a Gaussian with covariance $\bm{S}$ is strictly contained within the ellipsoid associated with $\bm{\Sigma}_{\bm{y}}$. 
Notably, the structure of the multivariate \gls*{GMM} allows for recursive decomposition along different directions, further enhancing its flexibility and applicability in high-dimensional settings.

\subsubsection{Low-order automatic domain splitting for automatic \gls*{GMM} decomposition}
Selecting an appropriate uncertainty propagation method is critical for effective stochastic optimisation.
Linear covariance propagation is a method used to propagate a Gaussian distribution $\bm{y}\sim \mathcal{N}\left(\Bar{\bm{y}}, \bm{\Sigma}_{\bm{y}}\right)$ of size $d$ through a linear function $\bm{h}$. The goal is to retrieve $\Bar{\bm{h}}$ and $\bm{\Sigma}_{\bm{h}}$ such that $\bm{h}\left(\bm{y}\right)\sim \mathcal{N}\left(\Bar{\bm{h}}, \bm{\Sigma}_{\bm{h}}\right)$. If $\nabla\bm{h}\left(\Bar{\bm{y}}\right)$ is the gradient of $\bm{h}$ at $\Bar{\bm{y}}$, then: $\Bar{\bm{h}}=\bm{h}\left(\Bar{\bm{y}}\right)$ and $\bm{\Sigma}_{\bm{h}}=\nabla\bm{h}\left(\Bar{\bm{y}}\right)\bm{\Sigma}_{\bm{y}}\nabla\bm{h}\left(\Bar{\bm{y}}\right)^{\textrm{T}}$. However, the linearity hypothesis is strong and often not satisfied, especially in the \gls*{CR3BP}. Other methods exist to do so, such as the unscented transform \citep{JulierUhlmann_1997_NEotKFtNS}, yet it also has limited resilience to high non-linearities.

\citet{LosaccoEtAl_2024_LOADSAfNUM} developed a \gls*{LOADS}-\gls*{GMM} method based on \gls*{DA} tools to automatically account for non-linearities as opposed to the linear covariance uncertainty propagation method.
To quantify non-linearity, they define the \gls*{NLI} $\eta\left(\mathcal{P}_{\bm{h}}\right)$, computed from the first- and second-order coefficients of the Taylor expansion $\mathcal{P}_{\bm{h}} = \bm{h}(\mathcal{P}_{\bm{y}})$. The polynomial $\mathcal{P}_{\bm{y}}$ is defined as $\mathcal{P}_{\bm{y}} = \Bar{\bm{y}}+\bm{L}_{\bm{y}}\delta \bm{y}$, where the Cholesky \citep{Cholesky_1910_SLRNDSDL} factorisation $\bm{L}_{\bm{y}}$ of $\bm{\Sigma}_{\bm{y}}$ serves as scale to represent the distribution $\bm{y}$, see \citet{LosaccoEtAl_2024_LOADSAfNUM} for the complete computation procedure.
The closer $\eta$ is to \num{0}, the less non-linear it is. Next, they define a threshold $\varepsilon_{\textrm{NLI}}$ such that if $\eta\leq\varepsilon_{\textrm{NLI}}$, $\bm{h}$ is considered linear, if not $\bm{h}$ is considered non-linear.

The \gls*{NLI} allows us to implement an automatic \gls*{GMM} decomposition algorithm. 
Once $\mathcal{P}_{\bm{h}}$ is computed, the \gls*{NLI} is evaluated and compared to a threshold $\varepsilon_{\textrm{NLI}}$. If $\eta \leq \varepsilon_{\textrm{NLI}}$, the function $\bm{h}$ is considered sufficiently linear, and linear covariance propagation is applied. Otherwise, the input distribution $\bm{y}$ is recursively decomposed into a \gls*{GMM} of size $M = 3$, following the procedure described in \Eq{eq:mulivariate_GMM}.
The decomposition direction is chosen based of the direction with higher directional \gls*{NLI} \citep{LosaccoEtAl_2024_LOADSAfNUM}. The distribution $\bm{y}$ can be recursively decomposed into a \gls*{GMM} of size $M$, denoted $\bm{y}\sim\textrm{GMM}\left(\bm{\alpha},\Bar{\bm{\mu}},\bm{S}\right)$.
Each mixand $\bm{y}_j \sim \mathcal{N}(\Bar{\bm{\mu}}_j, \bm{S}_j)$ is propagated through $\bm{h}$, yielding a Gaussian approximation $\bm{h}(\bm{y}_j) \sim \mathcal{N}(\Bar{\bm{h}}_j, \bm{\Sigma}_{\bm{h},j})$ with $\Bar{\bm{h}}_j=\bm{h}\left(\Bar{\bm{\mu}}_j\right)$ and $\bm{\Sigma}_{\bm{h},j}=\nabla\bm{h}\left(\Bar{\bm{\mu}}_j\right)\bm{S}_{j}\nabla\bm{h}\left(\Bar{\bm{\mu}}_j\right)^{\textrm{T}}$. The overall output distribution is then approximated as a \gls*{GMM}: $\bm{h}\left(\bm{y}\right)\sim\textrm{GMM}\left(\bm{\alpha},\left(\Bar{\bm{h}}_1,\cdots,\Bar{\bm{h}}_M\right),\left(\bm{\Sigma}_{\bm{h},1},\cdots,\bm{\Sigma}_{\bm{h},M}\right)\right)$.

\subsection{Chance-constrained optimal control problem formulation and transcription methods}
This section presents the formulation of chance constraints and the methods employed to address them.
Since any Gaussian random variable $\bm{y} \in \mathbb{R}^d$ satisfies $\mathbb{P}\left(\bm{y} \preceq \bm{0}\right) \in ]0, 1[$, a constraint of the form $\bm{y} \preceq \bm{0}$ cannot be enforced with absolute certainty ($\mathbb{P} = 1$). For this reason, chance constraints are used to handle such uncertainties. An acceptable failure rate $\beta \in ]0, 1[$ is defined, and the original constraints are reformulated as: $\mathbb{P}\left(\bm{y} \preceq \bm{0}\right) \geq 1-\beta$, which ensures that the constraint is satisfied with probability at least $1-\beta$.  
In the remainder of this paper, we define a failure as the event $\bm{y} \not\preceq \bm{0}$, \ie one of the components of $\bm{y}$ is positive, and refer to the associated failure risk as $\beta_{\textrm{R}} = 1 - \mathbb{P}\left(\bm{y} \preceq \bm{0}\right)$.

\subsubsection{Chance-constrained optimal control problem}
We now formulate the chance-constrained trajectory optimisation problem in the \gls*{DDP} framework. A similar formulation can be found in \citet{OzakiEtAl_2018_SDDPwUTfLTTD}. The states $\bm{x}_k$ and controls $\bm{u}_k$ are modeled as Gaussian variables: $\bm{x}_k \sim \mathcal{N}\left(\Bar{\bm{x}}_k, \bm{\Sigma}_{\bm{x},k}\right)$, $\bm{u}_k \sim \mathcal{N}\left(\Bar{\bm{u}}_k, \bm{\Sigma}_{\bm{u},k}\right)$.
The control policy consists of the nominal controls $\Bar{\bm{U}} = (\Bar{\bm{u}}_0, \cdots, \Bar{\bm{u}}_{N-1})$ and the associated linear feedback gains $\bm{K} = (\bm{K}_0, \cdots, \bm{K}_{N-1})$.  
The objective is to compute a policy $(\Bar{\bm{U}}, \bm{K})$ such that the failure risk is below $\beta$ and the $(1 - \beta)$-quantile of the cost $J$, denoted $\Tilde{J}$, is minimised. That is, the cost is below $\Tilde{J}$ with probability at least $1 - \beta$. This approach differs from minimising the expected cost, as done in \citet{OzakiEtAl_2020_TSOCfNCTOP}, for example. However, minimising the expected cost does not guarantee that the cost remains manageable in most cases, since deviations from the average cost are not controlled.

In the rest of this work, a subscript $k$ indicates evaluation at time step $k$, \eg $\bm{f}_k = \bm{f}\left(\bm{x}_k, \bm{u}_k\right)$, $\bm{g}_k = \bm{g}\left(\bm{x}_k, \bm{u}_k\right)$, and $\bm{g}_N = \bm{g}\left(\bm{x}_N, \bm{x}_t\right)$. 
Based on this formulation, the chance-constrained trajectory optimisation problem can be stated as follows:
\begin{equation}
    \label{eq:stochastic_opt_problem}
    \begin{aligned}
        & \min_{\Bar{\bm{U}}, \bm{K}} \Tilde{J}\left(\Bar{\bm{U}}, \bm{K}, \bm{x}_0, \bm{x}_t, \beta\right), \\
        & \text{subject to:} \\
        & J\left(\Bar{\bm{U}}, \bm{K}, \bm{x}_0, \bm{x}_t\right) = \sum_{k=0}^{N-1} l\left(\bm{x}_k, \bm{u}_k\right) + \phi\left(\bm{x}_N, \bm{x}_t\right), \\
        & \bm{x}_{k+1} = \bm{f}_k + \bm{w}_{\bm{x}, k+1}, \quad \forall k \in \mathbb{N}_{N-1}, \\
        & \bm{u}_k = \Bar{\bm{u}}_k + \bm{K}_k(\bm{x}_k - \Bar{\bm{x}}_k), \quad \forall k \in \mathbb{N}_{N-1}, \\
        & \mathbb{P}\left[J\left(\Bar{\bm{U}}, \bm{K}, \bm{x}_0, \bm{x}_t\right) \leq \Tilde{J}\left(\Bar{\bm{U}}, \bm{K}, \bm{x}_0, \bm{x}_t, \beta\right)\right] \geq 1 - \beta, \\
        & \mathbb{P}\left[\bigcap_{k=0}^{N} \left\{ \bm{g}_k \preceq \bm{0} \right\} \right] \geq 1 - \beta,
    \end{aligned}
\end{equation}
where the initial state $\bm{x}_0 \sim \mathcal{N}\left(\Bar{\bm{x}}_0, \bm{\Sigma}_{\bm{x},0}\right)$, the target state $\bm{x}_t \sim \mathcal{N}\left(\Bar{\bm{x}}_t, \bm{\Sigma}_{\bm{x},t}\right)$, and the process noise $\bm{w}_{\bm{x},k} \sim \mathcal{N}\left(\bm{0}, \bm{\mathcal{Q}}_{\bm{x},k}\right)$.

\subsubsection{Transcription methods}
To make the problem tractable, transcription methods are introduced to transform chance constraints into deterministic ones.  
Many stochastic optimisation techniques rely on such transcriptions to approximate chance constraints of the form $\mathbb{P}\left(\bm{y} \preceq \bm{0}\right) \geq 1 - \beta$ with tractable deterministic sufficient conditions $\mathcal{C}\left(\Bar{\bm{y}}, \bm{\Sigma}_{\bm{y}}, \beta\right)$. We then have: 
\begin{equation}
    \mathcal{C}\left(\Bar{\bm{y}}, \bm{\Sigma}_{\bm{y}}, \beta\right) \implies \mathbb{P}\left(\bm{y} \preceq \bm{0}\right) \geq 1 - \beta.
\end{equation}
Since $\mathcal{C}$ is only a sufficient condition, it introduces conservatism, i.e., additional safety margins to ensure constraint satisfaction. 

Transcription methods have been proposed \citep{BlackmoreEtAl_2011_CCOPPwO,RidderhofEtAl_2020_CCCCfLTMFTO,OguriLantoine_2022_SSCPfRLTTDuU,NakkaChung_2023_TOoCCNSSfMPuU}, although they often rely on too conservative bondings that may degrade performance \citep{BenedikterEtAl_2022_CAtCCwAtSLTTO,MarmoZavoli_2024_CCMfCCoLTIM} or they only tackle specific types of constraints. 
In addition, these transcription approaches are unable to generalise to multi-dimensional constraints. 
\citet{CalebEtAl_2025_CCTaFREfSTO} introduces a general multi-dimensional transcription method for Gaussian chance constraints with reduced conservatism and linear complexity.
If $\bm{y}$ is of size $d$, the first-order transcription is written as: 
\begin{equation}
    \label{eq:first_order_theorem}
    \Bar{\bm{y}} + \Psi_d^{-1}\left(\beta\right)\bm{\sigma}_{y} \leq \bm{0} 
    \implies \mathbb{P}\left(\bm{y}\preceq \bm{0}\right) \geq 1-\beta
\end{equation}
where $\bm{\sigma}_{y}$ is the vector of the positive square roots of the diagonal coefficients of $\bm{\Sigma}_{y}$ and $\Psi^{-1}_d\left(\beta\right) = \sqrt{\Phi^{-1}_d\left(1-\beta\right)}$, with $\Phi^{-1}_d$ the inverse \gls*{CDF} of chi-squared distribution with $d$ degrees of freedom \citep{AbramowitzStegun_1964_HoMFwFGaMT}.

\subsubsection{Risk estimation methods and conservatism metric}
\label{sec:risk_estimation_conservatism}
While deterministic transcriptions provide conservative guarantees, it is often desirable to estimate the actual failure risk $\beta_\mathrm{R} = 1 - \mathbb{P}(\bm{y} \preceq \bm{0})$ more accurately. 
Although the risk can, in principle, be estimated using numerical integration or Monte Carlo methods \citep{RobertCasella_2004_MCSM}, these approaches are generally intractable during the optimisation process.

\citet{CalebEtAl_2025_CCTaFREfSTO} introduce the so-called $d$\textendash{th}-order, a low-conservatism risk estimate.
Let us define $\Tilde{r}_1 \leq \cdots \leq \Tilde{r}_d$, the ordered components of the vector $\bm{r}=\left(-\frac{\bar{y}_i}{\sigma_{\bm{y},i}}\right)_{i\in\mathbb{N}_d^*}$.
If $\Bar{\bm{y}} \preceq \bm{0}$, then:
\begin{equation}
    \label{eq:def_beta_d}
    \beta_\mathrm{T} = 1 - \sum_{i=1}^d \left[\Psi_d(\Tilde{r}_i) - \Psi_d(\Tilde{r}_{i-1})\right] \max\left[0, 1 - \frac{1}{2} \sum_{j=1}^{i-1} I\left(\frac{d-1}{2}, \frac{1}{2}; \frac{\Tilde{r}_j}{\Tilde{r}_i}\right)\right] \geq \beta_\mathrm{R},
\end{equation}
where for $r\geq0$: $\Psi_d\left(R\right) = 1 - \Phi_d\left(R^2\right)$ and $I(a,b;x)$ denotes the regularised incomplete beta function: $I\left(a, b;x\right) = \frac{1}{B\left(a, b\right)} \int_{0}^{x} {t^{a-1}(1-t)^{b-1}dt}$, with $B(a,b)$ the Euler beta function \citep{AbramowitzStegun_1964_HoMFwFGaMT}.
Evaluating the risk estimation $\beta_{\rm T}$ requires $\mathcal{O}(d^2)$ operations and provides an upper bound of the failure risk $\beta_{\rm R}$ robust to high values of $d$. Indeed, if the conservatism of a risk estimation is defined as:
\begin{equation}
    \gamma(\beta_\mathrm{T}) = \frac{\beta_\mathrm{T}}{\beta_\mathrm{R}}  \sqrt{\frac{1 - \beta_\mathrm{R}^2}{1 - \beta_\mathrm{T}^2}},
\end{equation}
then the conservatism of the $d$\textendash{th} order estimator is stable and lower than \num{10}. See \citet{CalebEtAl_2025_CCTaFREfSTO} for additional precisions.

%% file: sections/3_l-soda.tex
\section{Linear stochastic optimisation with differential algebra}
\label{sec:LSODA}
This section presents a \gls*{DA}-based stochastic optimisation framework under the assumption that uncertainties are Gaussian and that their covariances are sufficiently small to justify local linear approximations of the system dynamics and constraints around their nominal values. 
After relaxing and transcribing the stochastic optimal control problem in \Eq{eq:stochastic_opt_problem}, we introduce the \gls*{L-SODA} solver, which leverages stochastic constrained \gls*{DA}-based \gls*{DDP} and a stochastic Newton solver.

\subsection{Transcription of the stochastic optimal control problem}

\subsubsection{Problem relaxation}
The chance constraints and cost function in \Eq{eq:stochastic_opt_problem} involve Gaussian variables across multiple time steps, making direct transcription challenging and potentially overly conservative.
To improve convergence of the \gls*{AUL} solver, the original problem is relaxed by decoupling the joint chance constraints into individual ones.
Therefore, this part of the problem:
\begin{equation}
    \begin{aligned}
        \mathbb{P}\left[J\left(\Bar{\bm{U}}, \bm{K}, \bm{x}_0, \bm{x}_t\right)  \leq \Tilde{J}\left(\Bar{\bm{U}}, \bm{K}, \bm{x}_0, \bm{x}_t, \beta\right)\right] & \geq 1 - \beta, \\
         \mathbb{P}\left[\bigcap_{k=0}^{N} \left\{ \bm{g}_k \preceq \bm{0} \right\} \right] &\geq 1 - \beta,
    \end{aligned}
\end{equation}
is relaxed into:
\begin{equation}
    \label{eq:relaxed_constraints}
    \begin{aligned}
        \mathbb{P}\left(l_k \leq \Tilde{l}_k\right) \geq& 1 - \beta,\ \forall k\in\mathbb{N}_{N-1}, \\
        \mathbb{P}\left(\phi_N \leq \Tilde{\phi}_N\right) \geq& 1 - \beta, \\
        \mathbb{P}\left( \bm{g}_k \preceq \bm{0} \right) \geq& 1 - \beta,\ \forall k\in\mathbb{N}_{N}, \\
    \end{aligned}
\end{equation}
where $\Tilde{l}_k$ and $\Tilde{\phi}_N$ are the $1-\beta$ quantiles of $l_k$ and $\phi_N$. Note that this version of the problem is a necessary condition for the previous one. For example, if a portion $\beta$ of the trajectories violate the constraint at step $k$ and another disjoint $\beta$ violate it at step $k+1$, each individual constraint may appear satisfied, while the overall failure probability exceeds $2\beta$, rendering the trajectory infeasible.
This formulation is only used to provide a solution to then warm start the complete stochastic optimal control problem.

\subsubsection{Inequalities and costs}
Based on the relaxed formulation, we now describe how to transcribe the stochastic costs and constraints into tractable deterministic ones. 
The chance constraints or costs are functions of the state and the control\footnote{Those that are functions of the state and the target state can be transcribed using a similar procedure.} $\mathbb{P}\left(\bm{h}\left(\bm{x}_k, \bm{u}_k\right) \preceq \bm{0} \right) \geq 1-\beta$, where $\bm{h}$ can be considered linear, $\delta\bm{x}_k=\bm{x}_k-\Bar{\bm{x}}_k$, and $\delta\bm{u}_k=\bm{u}_k-\Bar{\bm{u}}_k$.
Since the control policy is affine, we substitute $\delta\bm{u}_k = \bm{K}_k \delta\bm{x}_k$, yielding a linearised form of $\bm{h}$ in terms of $\delta\bm{x}_k$:
$\bm{h}\left(\bm{x}_k, \bm{u}_k\right) = \bm{h}\left(\Bar{\bm{x}}_k, \Bar{\bm{u}}_k\right) + \left(\bm{h}_{\bm{x},k} + \bm{h}_{\bm{u},k}\bm{K}_k\right)\delta\bm{x}_k$.
The vector $\delta\bm{x}_k\sim\mathcal{N}\left(\bm{0},\bm{\Sigma}_{\bm{x},k}\right)$, therefore, $\bm{h}\left(\bm{x}_k, \bm{u}_k\right)=\bm{y}\sim\mathcal{N}\left(\Bar{\bm{y}},\bm{\Sigma}_{y}\right)$, where: $\Bar{\bm{y}} = \bm{h}\left(\Bar{\bm{x}}_k, \Bar{\bm{u}}_k\right)$, $ \bm{\Sigma}_{y} = \left(\bm{h}_{\bm{x},k} + \bm{h}_{\bm{u},k}\bm{K}_k\right)\bm{\Sigma}_{\bm{x},k}\left(\bm{h}_{\bm{x},k} + \bm{h}_{\bm{u},k}\bm{K}_k\right)^{\textrm{T}}$.
Therefore, all chance constraints of \Eq{eq:relaxed_constraints} can be considered Gaussian chance constraints, and using a transcription $\mathcal{T}$ such as the first-order transcription, we now have:
\begin{equation}
    \mathcal{T}\left(\bm{g}_k,\beta\right) \preceq \bm{0},\ \forall k\in\mathbb{N}_{N} \implies  \mathbb{P}\left( \bm{g}_k \preceq \bm{0} \right) \geq 1 - \beta,\ \forall k\in\mathbb{N}_{N}.
\end{equation}
Regarding the costs of \Eq{eq:relaxed_constraints}, we only need to replace each $l_k$ and $\phi_N$ by $\mathcal{T}\left(l_k,\beta\right)$ and $\mathcal{T}\left(\phi_N,\beta\right)$ in the \gls*{AUL} solver so that we have $\mathbb{P}\left(l_k \leq \Tilde{l}_k\right) \geq 1 - \beta$ and $\mathbb{P}\left(\phi_N \leq \Tilde{\phi}_N\right) \geq 1 - \beta$.

\subsubsection{Terminal constraints}
Enforcing equality constraints in a probabilistic setting, such as $\mathbb{P}(\bm{y} = \bm{0}) \geq 1 - \beta$, is infeasible for continuous distributions like Gaussian ones, since the probability of any exact value is zero.
However, terminal equality constraints are often required to ensure the final state $\bm{x}_N$ matches the target state $\bm{x}_t$. 
These constraints are often strong constraints while we only need to ensure the spacecraft is in a given area with a certain probability. 

Instead, we define terminal constraints in terms of a confidence region: the final state $\bm{x}_N$ must lie within a Mahalanobis distance from the target mean $\Bar{\bm{x}}_t$, corresponding to a confidence level $1 - \beta_t$ \citep{Mahalanobis_1936_OtGDiS}:
\begin{equation}
    \label{eq:terminal_constraints}
    g_N=\left(\bm{x}_N-\Bar{\bm{x}}_t\right)^{\textrm{T}}\bm{\Sigma}_t^{-1}\left(\bm{x}_N-\Bar{\bm{x}}_t\right) - \Phi_{N_x}^{-1}\left(1-\beta_t\right)\leq 0.
\end{equation}
This constraint ensures the distribution $\bm{x}_N$ is "within" that of $\bm{x}_t$ with a probability $\beta_t$ and it can also be transcribed as a function of $\bm{x}_N$, a Gaussian distribution.

\subsection{Linear stochastic \gls*{AUL} solver}
\subsubsection{Propagation and differentiation of stochastic quantities using \gls*{DA}}

To enable stochastic optimisation with \gls*{DDP}, both the propagation of uncertainties and the computation of derivatives must be performed efficiently. Leveraging the \gls*{DA} framework, we use Taylor expansions to compute the means and covariances of stochastic quantities, as well as their first- and second-order derivatives. Since the Taylor expansions of the dynamics $\bm{f}_k$, the stage cost $\hat{l}_k$, and the terminal cost $\hat{\phi}_N$ are already computed in the constrained \gls*{DA}-based \gls*{DDP} solver \citep{CalebEtAl_2025_APBCSfFOLTTO}, these quantities can be directly reused throughout the optimisation process.
For example, consider the system dynamics: $\bm{x}_{k+1} = \bm{f}(\bm{x}_k, \bm{u}_k) + \bm{w}_{\bm{x},k+1}$.
Using linear covariance propagation, the mean and covariance of $\bm{x}_{k+1}$ are given by: 
\begin{equation}
    \Bar{\bm{x}}_{k+1} = \overline{\mathcal{P}_{\bm{f},k}}, \quad \bm{\Sigma}_{\bm{x},k+1} = \left(\bm{f}_{\bm{x},k} + \bm{K}_k \bm{f}_{k,\bm{u}}\right) \bm{\Sigma}_{\bm{x},k} \left(\bm{f}_{\bm{x},k} + \bm{K}_k \bm{f}_{k,\bm{u}}\right)^{\textrm{T}} + \bm{\mathcal{Q}}_{\bm{x},k+1},
\end{equation}
This computation involves only matrix operations and is therefore highly efficient.

Performing \gls*{DDP} optimisation also requires the derivatives of the transcribed costs and constraints. These are obtained by expanding any stochastic function $\bm{h}(\bm{y})$ around its mean $\Bar{\bm{y}}$ using a Taylor expansion: 
\begin{equation}
    \mathcal{P}_{h} = h(\Bar{\bm{y}}) + \nabla \bm{h}(\Bar{\bm{y}}) \delta \bm{y} + \frac{1}{2} \delta \bm{y}^{\textrm{T}} \frac{\partial^2 \bm{h}}{\partial \bm{y}^2}(\Bar{\bm{y}}) \delta \bm{y} + \dotsc,
\end{equation}
here with $h$ scalar for the sake of conciseness.
Differentiating this expansion yields: $\mathcal{P}_{\nabla \bm{h}} = \nabla \bm{h}(\Bar{\bm{y}}) + \frac{\partial^2 \bm{h}}{\partial \bm{y}^2}(\Bar{\bm{y}}) \delta \bm{y} + \dotsc$.
The covariance of $\bm{h}(\bm{y})$ can then be expressed as: $\mathcal{P}_{\bm{\Sigma}_{\bm{h}}} = \mathcal{P}_{\nabla \bm{h}} \bm{\Sigma}_{\bm{y}} \mathcal{P}_{\nabla \bm{h}}^{\textrm{T}}$.
Since the transcriptions $\mathcal{T}(\bm{h}, \beta)$ depend on the mean and covariance of $\bm{h}$, they can also be expanded as Taylor polynomials: $\mathcal{P}_{\mathcal{T}(\bm{h})} = \mathcal{T}(\mathcal{P}_{\bm{h}}, \beta)$.
Because \gls*{DDP} requires second-order derivatives and the covariance expansion involves differentiating the Jacobian, the Taylor expansions must be computed to at least third order. This is in contrast to the \gls*{DADDy} solver, which only requires second-order expansions \citep{CalebEtAl_2025_APBCSfFOLTTO}. Although this increases the computational cost, it also improves robustness by enlarging the convergence radius of the expansions. For the remainder of this work, all expansions are performed to order three.

\subsubsection{Implementation of the solver}
This subsection details the practical implementation of the \gls*{L-SODA} solver. From this point onward, the state $\bm{x}_k$ and control $\bm{u}_k$ are treated as distributions, represented respectively as $\bm{x}_k = (\Bar{\bm{x}}_k, \bm{\Sigma}_{\bm{x},k})$ and $\bm{u}_k = (\Bar{\bm{u}}_k, \bm{K}_k)$. Accordingly, $\bm{X}$ and $\bm{U}$ denote the sequences of nominal states with covariances and nominal controls with feedback gains.
\Alg{alg:AUL_linear_stochastic} presents the linear stochastic \gls*{AUL} solver.
\begin{algorithm}[h]
    \caption{Linear stochastic \gls*{AUL} solver}
    \label{alg:AUL_linear_stochastic}
    \begin{algorithmic}
        \State \textbf{Input:} $\varepsilon_{\textrm{DDP}}$, $\varepsilon_{\textrm{AUL}}$, $\beta$, $\bm{x}_0$, $\bm{x}_t$, $\bm{U}_0$, $\bm{f}$, $l$, ${\phi}$, $\bm{g}$, $\bm{g}_t$, $\mathcal{T}$
        \State $\bm{U} \gets \bm{U}_0$
        \State Initialise penalties and dual states
        \State Initialise augmented costs functions $\hat{l}$ and $\hat{\phi}$
        \State Compute $\bm{X}$, $\mathcal{P}_{\bm{f}}$, $\mathcal{P}_{\mathcal{T}(\hat{l})}$, and $\mathcal{P}_{\mathcal{T}(\hat{\phi})}$, $\bm{G}= (\bm{g}_0, \cdots, \bm{g}_{N})$
        \State $g_{\max}\gets\max \bm{G}$
        \While{$g_{\max}>\varepsilon_{\textrm{AUL}}$}
            
            \State $\Tilde{J}\gets -\infty$
            \State Compute $\Tilde{J}^*$
            \While{$\Tilde{J}^* < \Tilde{J}\land |\Tilde{J}-\Tilde{J}^*|>\varepsilon_{\textrm{DDP}}$}
                \State $\Tilde{J}\gets \Tilde{J}^*$
                \State $\bm{A}, \ \bm{K} \gets$ BackwardSweep$\left(\bm{X}, \bm{U}, \bm{x}_t, \mathcal{P}_{\bm{f}},  \mathcal{P}_{\mathcal{T}(\hat{l})}, \mathcal{P}_{\mathcal{T}(\hat{\phi})}\right)$                
                \State $\Tilde{J}^*,\ \bm{X}^*,\ \bm{U}^*,\ \mathcal{P}_{\bm{f}}, \ \mathcal{P}_{\mathcal{T}(\hat{l})}, \ \mathcal{P}_{\mathcal{T}(\hat{\phi})} \gets$ ForwardPass$\left(\bm{X}, \bm{U}, \bm{x}_t, \bm{A},\bm{K}, \mathcal{P}_{\bm{f}}, \mathcal{T}(\hat{l}), \mathcal{T}(\hat{\phi})\right)$    
                \State $\bm{X},\ \bm{U} \gets \bm{X}^*,\ \bm{U}^*$     
            \EndWhile
            \State $\Tilde{J}\gets \Tilde{J}^*$
            \State Compute $\bm{G}$
            \State $g_{\max}\gets \max\bm{G}$
            \State Update penalties and dual states
            \State Update augmented costs functions $\hat{l}$ and $\hat{\phi}$
        \EndWhile
        \State \textbf{Return:}  $\Tilde{J}$, $\bm{X}$, $\bm{U}$, $g_{\max}$
    \end{algorithmic}
\end{algorithm}
The BackwardSweep function is the one from \citet{CalebEtAl_2025_APBCSfFOLTTO} and the ForwardPass function is presented in \Alg{alg:forward_pass_DA_stochastic}. This function propagates both the nominal trajectory and the associated uncertainties.
\begin{algorithm}[h]
    \caption{Forward pass}
    \label{alg:forward_pass_DA_stochastic}
    \begin{algorithmic}
        \State \textbf{Input:} $\varepsilon_{\textrm{DA}}$, $\bm{X}$, $\bm{U}$, $\bm{x}_t$, $\bm{A}= (\bm{a}_0, \cdots, \bm{a}_{N-1})$, $\bm{K}$, $\mathcal{P}_{\bm{f}}$, $l$, $\phi$
        \State $k,\ J^*,\  \bm{x}_k^* \gets 0, \ 0, \ \bm{x}_0$
        \While{$k\leq N-1$}
            \State $\delta\Bar{\bm{x}}_{k}^*\gets\Bar{\bm{x}}_{k}^*-\Bar{\bm{x}}_{k}$ 
            \State $\delta\Bar{\bm{u}}_k^* \gets \bm{a}_k + \bm{K}_k\delta\Bar{\bm{x}}_k^{*}$ 
            \State $\Bar{\bm{u}}_k^* \gets\Bar{\bm{u}}_k + \delta\Bar{\bm{u}}_k^* $ 
            \State Compute $R_{\varepsilon_{\textrm{DA}},k}$
            \If{$\sqrt{\|\delta\bm{x}_k^*\|^2 + \|\delta\bm{u}_k^*\|^2}<R_{\varepsilon_{\textrm{DA}},k}$}
                \State $\mathcal{P}_{f,k} \gets \mathcal{P}_{f,k}\left(\delta\Bar{\bm{x}}_k^* + \delta\bm{x}, \delta\Bar{\bm{u}}_k^* + \delta\bm{u}\right)$
            \Else
                \State $\mathcal{P}_{f,k} \gets \bm{f}\left(\Bar{\bm{x}}_k^* + \delta\bm{x}, \Bar{\bm{u}}_k^*+ \delta\bm{u}\right)$
            \EndIf
            \State $\mathcal{P}_{l,k}\gets l\left(\Bar{\bm{x}}_k^* + \delta\bm{x}, \Bar{\bm{u}}_k^* + \delta\bm{u}\right)$
            \State $\Bar{\bm{x}}_{k+1}^*,\ J^* \gets \overline{\mathcal{P}_{f,k}},\ J^* + \overline{\mathcal{P}_{l,k}}$ 
            \State Compute $\bm{\Sigma}_{\bm{x},k+1}^{*}$
            \State $k \gets k + 1$
        \EndWhile
        \State $\mathcal{P}_{\phi}\gets \phi\left(\Bar{\bm{x}}_N^* + \delta\bm{x}, \bm{x}_t\right)$
        \State $J^*\gets J^* + \overline{\mathcal{P}_{\phi}}$
        \State \textbf{Return:} $J^*$, $\bm{X}^*$, $\bm{U}^*$, $\mathcal{P}_{\bm{f}}$, $\mathcal{P}_{l}$, $\mathcal{P}_{\phi}$
    \end{algorithmic}
\end{algorithm}
Throughout this implementation, \gls*{DA} expansions serve three key purposes: automatic differentiation, dynamics approximation, and uncertainty propagation. This unified use of Taylor models significantly reduces computational overhead while maintaining accuracy.

\subsection{Stochastic adaptive Newton solver}
Once a solution $\bm{U}$ of the relaxed problem is computed using the stochastic \gls*{AUL} solver, a solution polishing method is used to make it a solution of the initial chance-constrained problem.
This section presents the implementation of a stochastic adaptive Newton solver.

\subsubsection{Adaptive first-order method}
The main idea behind this solver is that the vector of constraints $\bm{\Gamma}$ is of high-dimension and transcribing it using spectral radius of first-order transcription would result in important conservatism, thus, sub-optimality. 
In this part, $\bm{\Gamma}$ refers only to the stochastic constraints, \ie all except the continuity constraints. 
Since the chance constraints are of high dimensions, the first-order transcription alone as a large conservatism. Therefore, this solver uses the $d$-th order risk estimation $\beta_{\textrm{T}}$ combined with the first-order transcription to satisfy constraints in the least conservative way.
To do so, note that $\bm{\Gamma}\left(\bm{X},\bm{U}\right)$ is a function of Gaussian variables with small covariance, therefore, it can be considered a Gaussian variable: $\bm{\Gamma}\left(\bm{X},\bm{U}\right)\sim\mathcal{N}\left(\Bar{\bm{\Gamma}},\bm{\Sigma}_{\bm{\Gamma}}\right)$. This vector can be transcribed with the first-order transcription: $\mathcal{T}\left(\bm{\Gamma},\beta\right)=\Bar{\bm{\Gamma}}+\Psi_{N_{\Gamma}}^{-1}\left(\beta\right)\bm{\sigma}_{\bm{\Gamma}}$. The adaptive Newton method consist in modifying the transcription into:
\begin{equation}
    \mathcal{T}_{\tau}\left(\bm{\Gamma},\beta\right)=\Bar{\bm{\Gamma}}+\tau\Psi_{N_{\Gamma}}^{-1}\left(\beta\right)\bm{\sigma}_{\bm{\Gamma}},
\end{equation}
where $\tau\in]0,1]$. The value $\tau$ depends on the estimated failure risk $\beta_{\textrm{T}}$ and the target risk $\beta$. At each iteration of the Newton solver, $\tau$ increases of $s=\frac{(\beta_{\textrm{T}}-\beta)^2}{2}$ if $\beta_{\textrm{T}}>2\beta$, \ie constraints are not conservative enough, or decreases of $s$ if $2\beta_{\textrm{T}}<\beta$, the constraints are too conservative. This specific value of $s$ was empirically found to be a good compromise to make small adjustments when $\beta_{\textrm{T}}\approx\beta$ and large ones when they are different, typically \SI{100}{\percent} and \SI{5}{\percent}.

\subsubsection{Computing the $d$-th order risk estimation}
The $d$-th order risk estimation of a distribution $\bm{y}$ is computed thanks to the square root of the diagonal of the covariance $\bm{\sigma}_{\bm{y}}$ and the mean $\Bar{\bm{y}}$. In the case of $\bm{\Gamma}$, of gradient $\bm{\Delta}$, retrieving $\Bar{\bm{\Gamma}}$ is straightforward. However, the naive computation of $\bm{\Sigma}_{\bm{\Gamma}}=\bm{\Delta}\bm{C}\bm{\Delta}^{\textrm{T}}$, where $\bm{C}$ is a block diagonal matrix composed by $\bm{\Sigma}_{\bm{x},0},\bm{\Sigma}_{\bm{u},0}\dotsc,\bm{\Sigma}_{\bm{x},N}$, is of complexity $\mathcal{O}(N_{\Gamma}N_{Y}(N_{\Gamma}+ N_{Y}))$ (assuming naive matrix multiplication complexity). Moreover, $\bm{\Delta}$ is a also a block diagonal matrix containing the gradients of each chance-constraint $\bm{g}_k$ \citep{CalebEtAl_2025_APBCSfFOLTTO}. Thus, this procedure can be accelerated. Indeed, the computation of $\bm{\sigma}_{\bm{\Gamma}}$ can be reduced to computing the gradient of each $\bm{g}_k$ denoted $\bm{\Delta}_k=\bm{g}_{\bm{x},k}+\bm{K}_k\bm{g}_{\bm{u},k}$. Then, computing each product $\bm{\Pi}_k=\bm{\Sigma}_{\bm{x},k}\bm{\Delta}_k^{\textrm{T}}$ is an operation of complexity $\mathcal{O}(N_x^2N_{g})$. 
Finally, the $i$-th diagonal term of $\bm{\Sigma}_{\bm{\Gamma}}$ is retrieved by finding the correct pair $(q,r)$ such that $i=qN_g + r$ and multiplying the $r$-th column of $\bm{\Pi}_q$ with the $r$-th row of $\bm{\Delta}_q$.
This operation is of complexity $\mathcal{O}(N_x)$, therefore, computing the full vector $\bm{\sigma}_{\bm{\Gamma}}$ has a complexity $\mathcal{O}(NN_x^2N_g)$.
For typical values: $N_x=7$, $N_u=3$, $N_g=3$, and $N=100$, the naive computation of $\bm{\Sigma}_{\bm{\Gamma}}$ requires on the order of \num{4e8} operations, compared to only \num{1.5e4} when exploiting sparsity. This dramatic reduction highlights the importance of using sparse matrix structures for the repeated evaluation of $\beta_{\textrm{T}}$.
\Alg{alg:linear_pn_method_stochastic} presents the linear stochastic Newton solver.
\begin{algorithm}[h]
    \caption{Linear stochastic Newton solver}
    \label{alg:linear_pn_method_stochastic}
    \begin{algorithmic}
        \State \textbf{Input:} $\beta$, $\varepsilon_{\textrm{N}}$, $\bm{X}_0$, $\bm{U}_0$, $\bm{x}_t$, $\bm{f}$, $l$, ${\phi}$, $\bm{g}$, $\bm{g}_t$
        \State Compute $\bm{Y}$
        \State $d_{\max},\ \zeta,\ \tau \gets +\infty,\ 0.5,\ \tau_{\textrm{init}}$   
        \While{$d_{\max}>\varepsilon_{\textrm{N}}$}
            \State $\mathcal{P}_{\mathcal{T}\left(\bm{\Gamma}\right)}\gets \mathcal{T}_{\tau}\left(\bm{\Gamma} \left(\bm{Y} + \delta\bm{Y}\right),\beta\right)
             $
            \State $\bm{d}\gets\overline{\mathcal{P}_{\mathcal{T}\left(\bm{\Gamma}\right)}}$
            \State Compute $\bm{\Delta}$ and $\bm{\Delta}\bm{\Delta}^{\textrm{T}}$ from $\mathcal{P}_{\mathcal{T}\left(\bm{\Gamma}\right)}$
            \State $\bm{\Pi}\gets$BlockCholesky$\left(\bm{\Delta}\bm{\Delta}^{\textrm{T}}\right)$ 
            \State $r\gets +\infty$ 
            \While{$d_{\max}>\varepsilon_{\textrm{N}}\land r>\varepsilon_{\textrm{CV}}$}
                \State $\delta\bm{Y}^*\gets -\bm{\Delta}^{{\textrm{T}}}{\left(\bm{\Pi}^{{\textrm{T}}}\right)^{-1}\bm{\Pi}^{-1} \bm{d}} $
                \State $d_{\max}^*,\ \kappa \gets +\infty,\ 1$ 
                \While{$d_{\max}^*>d_{\max}$}
                    \State $\bm{d}^*\gets\mathcal{P}_{\mathcal{T}\left(\bm{\Gamma}\right)}\left(\kappa\delta\bm{Y}^*\right)$
                    \State $d_{\max}^*,\ \kappa\gets \max |\bm{d}^*|, \ \zeta\kappa$
                \EndWhile
                \State $\bm{Y},\ \bm{d},\ \mathcal{P}_{\mathcal{T}\left(\bm{\Gamma}\right)}\gets\bm{Y} + \frac{\kappa}{\zeta}\delta\bm{Y}^*,\  \bm{d}^*,\ \mathcal{P}_{\mathcal{T}\left(\bm{\Gamma}\right)}\left(\frac{\kappa}{\zeta}\delta\bm{Y}^* + \delta\bm{Y}\right)$
                \State $r,\ d_{\max}\gets \log d_{\max}^*/\log d_{\max},\ \max |\bm{d}|$
            \EndWhile
            \State $\beta_{\textrm{T}}\gets \text{RiskEstimation}(\Bar{\bm{Y}},\bm{\Sigma})$
            \State $\tau\gets \text{UpdateTau}(\tau,\beta,\beta_{\textrm{T}})$
        \EndWhile
        \State \textbf{Return:}  $\bm{X}$, $\bm{U}$, $d_{\max}$, $\beta_{\textrm{T}}$
    \end{algorithmic}
\end{algorithm}

%% file: sections/4_soda.tex
\section{Stochastic optimisation with differential algebra}
\label{sec:soda}
Now that stochastic optimisation problems can be addressed under linear or weakly non-linear assumptions, this section introduces methods for solving non-linear stochastic optimal control problems in the fully non-linear setting using \gls*{LOADS}-based \gls*{GMM}. 
We first present theoretical results and introduce a sequential stochastic optimisation method based on a novel failure risk allocation strategy. The remainder of this section then details how the non-linear problem is decomposed into a set of sub-problems.

\subsection{Failure risk allocation method}
\label{sec:soda_fam}
In this section, we consider a random variable $\bm{y}\sim\textrm{GMM}\left(\bm{\alpha},\Bar{\bm{\mu}},\bm{S}\right)$ of $M$ mixands.
The objective is to satisfy the non-Gaussian chance constraint: $\mathbb{P}\left(\bm{y}\preceq\bm{0}\right)\geq 1-\beta$. For each Gaussian mixand $\bm{y}_j \sim \mathcal{N}(\Bar{\bm{\mu}}_j, \bm{S}_j)$, we define $\beta_{\textrm{T},j}$ a conservative (tractable) estimation of the failure risk, $\beta_{\textrm{R},j}$ the true (intractable) failure risk, and $\beta_{\textrm{R}}$ the global failure risk of the \gls*{GMM}.
First, observe that: $\beta_{\textrm{R}} = \sum_{j=1}^{M}{\alpha_j\beta_{\textrm{R},j}}$, from this result, we observe that if each mixand is such that $\beta_{\textrm{T},j}\leq \beta$ and if we define $\beta_{\textrm{T}} = \sum_{j=1}^{M}{\alpha_j\beta_{\textrm{T},j}}$, then $\beta_{\textrm{T}}\leq\sum_{j=1}^{M}{\alpha_j\beta}=\beta$. Moreover, since by design $\beta_{\textrm{R},j}\leq\beta_{\textrm{T},j}$, we obtain $\beta_{\textrm{R}}\leq\beta$. Therefore, ensuring that each mixand individually satisfies its chance constraints guarantees that the entire \gls*{GMM} satisfies the global chance constraints.

However, in sequential optimisation, the target failure risk $ \beta $ can be dynamically relaxed to reduce conservatism. For example, if all $ \alpha_j = 1/M $ and the first mixand satisfies $ \beta_{\textrm{T},1} = 0 $, then the constraint on the second mixand can be relaxed to $ \mathbb{P}(\bm{y}_2 \preceq \bm{0}) \geq 1 - 2\beta $, while still ensuring $ \beta_{\textrm{R}} \leq \beta $. In other words, unused failure risk margins from previous mixands can be reallocated to subsequent ones, facilitating a more efficient optimisation process.
Let us define $\delta_j = \beta^*_j-\beta_{\textrm{T},j}$ and $\Delta_j = \beta^*_j-\beta_{\textrm{R},j}$ for $j\in\mathbb{N}_M^*$, $\beta^*_{j+1} = \beta + \dfrac{\alpha_j}{\alpha_{j+1}}\delta_j$ for $j\in\mathbb{N}_{M-1}^*$, and $\beta^*_1=\beta$. 
Note that $\beta^*_j\geq\beta$ and that $\Delta_j$ is not easily tractable in practice as opposed to $\delta_j$ and $\beta^*_j$ that can be accessed using failure risk estimation.
We propose a result to sequentially relax chance constraints while remaining robust with a target probability $\beta$.
\begin{theorem}
    \label{thm:FAM}
    $\forall j\in\mathbb{N}_M^*,\ \beta_{\textrm{T},j}\leq\beta^*_j\implies \beta_{\textrm{R}}\leq\beta$.
\end{theorem}
\begin{proof}
    By hypothesis $\delta_j= \beta^*_j-\beta_{\textrm{T},k}\geq 0$, then, by design of a risk estimation method $\beta_{\textrm{T},j}-\beta_{\textrm{R},j}\geq 0$. Thus: $\Delta_j = \beta^*_j-\beta_{\textrm{R},j} = \delta_j +\beta_{\textrm{T},j}-\beta_{\textrm{R},j} \geq \delta_j\geq 0$.
    Therefore:
    \begin{equation}
        \begin{aligned}
            \beta - \beta_\textrm{R} =&\sum_{j=1}^{M}{\alpha_j\left(\beta-\beta_{\textrm{R},j}\right)} = \alpha_1\left(\beta_1^{*}-\beta_{\textrm{R},1}\right) + \sum_{j=2}^{M}{\alpha_j\left(\beta^*_{j} - \dfrac{\alpha_{j-1}}{\alpha_{j}}\delta_{j-1}-\beta_{\textrm{R},j}\right)}  \\ = & \alpha_1\Delta_1 + \sum_{j=2}^{M}{\alpha_j\left(\Delta_{j} - \dfrac{\alpha_{j-1}}{\alpha_{j}}\delta_{j-1}\right)} \\
            = & \alpha_1\Delta_1 + \sum_{j=2}^{M}{\alpha_j\Delta_{j} } - \sum_{j=2}^{M}{\alpha_{j-1}\delta_{j-1}} = \sum_{j=1}^{M-1}{\alpha_j\Delta_{j} } - \sum_{j=1}^{M-1}{\alpha_{j}\delta_{j}} + \alpha_M\Delta_{M} \\
            = & \sum_{j=1}^{M-1}{\alpha_{j}\left(\Delta_{j}-\delta_{j}\right)} + \alpha_M\Delta_{M}\\
        \end{aligned}
    \end{equation}
    Finally, since $\Delta_j\geq\delta_j\geq 0$, we have $\beta - \beta_\textrm{R}\geq 0$.
\end{proof}
This result shows that if for each mixand $\beta_{\textrm{T},j}\leq\beta^*_j$, then, the chance constraints are satisfied for the \gls*{GMM} $\bm{y}$. The proof of \Thm{thm:FAM} also shows that the quantity $\Delta = \beta - \beta_{\textrm{R}}$ computed using the failure allocation method is smaller than $\Delta'=\sum_{j=1}^{M}{\alpha_j\Delta_j}$, the same difference computed by satisfying simply $\beta_{\textrm{T},j}\leq\beta^*_{j}=\beta$.

\subsection{\Gls*{GMM}-based problem decomposition and optimisation initialisation}
\subsubsection{\Gls*{GMM} decomposition of the trajectory optimisation problem}
We now focus on how to decompose a non-linear stochastic optimal control problem into sub-problems that can be efficiently handled by the linear solver. 
Consider a trajectory $\left(\bm{X}, \bm{U}\right)$ during the optimisation process. Before entering the \gls*{DDP} solver, the \gls*{NLI} is evaluated at each time step. If the maximum \gls*{NLI} exceeds a given threshold $\varepsilon_{\textrm{NLI}}$, the initial distribution $\bm{x}_0$ is decomposed into three mixands $\bm{x}_{0,-1}, \bm{x}_{0,0}, \bm{x}_{0,1}$ along the most non-linear dimension $j$, following the procedure of \citet{LosaccoEtAl_2024_LOADSAfNUM} shown in \Sec{sec:GMM}.
The trajectories corresponding to each mixand are generated using the \gls*{DA}-based expansion of the dynamics: $\Bar{\bm{u}}_{k,s} = \Bar{\bm{u}}_{k} + \bm{K}_k \left(\Bar{\bm{x}}_{k,s} - \Bar{\bm{x}}_{k} \right)$, $\mathcal{P}_{\bm{f},k,s} = \mathcal{P}_{\bm{f},k}\left(\Bar{\bm{x}}_{k,s} - \Bar{\bm{x}}_{k} + \delta \bm{x},\, \Bar{\bm{u}}_{k,s} - \Bar{\bm{u}}_{k} + \delta \bm{u} \right)$, $\Bar{\bm{x}}_{k+1,s} = \overline{\mathcal{P}_{\bm{f},k,s}}$, where $s \in \{-1, 0, 1\}$. The covariance of each distribution $\bm{x}_{k,s}$ is propagated using a linearised model, based on the gradient of $\mathcal{P}_{\bm{f},k,s}$. The gain $\bm{K}_k$ is reused for each $\bm{u}_{k,s}$.
This yields three distinct trajectories $\left(\bm{X}_s, \bm{U}_s\right)$. The lateral trajectories $s = -1$ and $s = 1$ are added to a list of trajectories $\mathcal{L}$ to be optimised later. The central trajectory $\left(\bm{X}_0, \bm{U}_0\right)$ is recursively split using the same criterion until its highest \gls*{NLI} falls below $\varepsilon_{\textrm{NLI}}$. At that point, it is optimised using the \gls*{DDP} solver under the weakly non-linear assumption, as described in \Sec{sec:LSODA}.
To avoid an excessive number of mixands, a lower bound $\alpha_{\min}$ is enforced on the \gls*{GMM} weights such that no mixand has weight $\alpha_j \leq \alpha_{\min}$. Once optimised, the trajectories are added to the set $\mathcal{L}^*$.

\subsubsection{First-guess generation for sub-problem optimisation}

To accelerate the optimisation of each trajectory in the list $\mathcal{L}$, we leverage the already optimised trajectories in $\mathcal{L}^*$ as first guesses.
Given a trajectory $\left(\bm{X}, \bm{U}\right)$ with initial condition $\bm{x}_0 \sim \mathcal{N}\left(\Bar{\bm{x}}_0, \bm{\Sigma}_{\bm{x},0}\right)$, we search in $\mathcal{L}^*$ for the trajectory with initial mean $\Bar{\bm{x}}_0^*$ that minimises the squared Mahalanobis distance: $\left(\Bar{\bm{x}}_0 - \Bar{\bm{x}}_0^*\right)^{\textrm{T}} \bm{\Sigma}_{\bm{x},0}^{-1} \left(\Bar{\bm{x}}_0 - \Bar{\bm{x}}_0^*\right)$.
The corresponding optimised trajectory $\left(\bm{X}^*, \bm{U}^*\right)$ is selected, and its control sequence is used to initialise the optimisation of the new trajectory as follows: $\Bar{\bm{u}}_{k} = \Bar{\bm{u}}_{k}^* + \bm{K}_k^* \left(\Bar{\bm{x}}_{k} - \Bar{\bm{x}}_{k}^* \right)$, $\mathcal{P}_{\bm{f},k} = \mathcal{P}_{\bm{f},k}^*\left(\Bar{\bm{x}}_{k} - \Bar{\bm{x}}_{k}^* + \delta \bm{x},\, \Bar{\bm{u}}_{k} - \Bar{\bm{u}}_{k}^* + \delta \bm{u} \right)$, $\Bar{\bm{x}}_{k+1} = \overline{\mathcal{P}_{\bm{f},k}}$, with $\bm{K}_k = \bm{K}_k^*$, and state covariances propagated using the dynamics $\mathcal{P}_{\bm{f},k}$.

\subsection{Overview of the stochastic optimisation with differential algebra solver}
\Alg{alg:AUL_stochastic} describes the stochastic \gls*{AUL} solver, which incorporates a \gls*{GMM}-based decomposition of trajectories exhibiting excessive \gls*{NLI}. At each iteration, a candidate trajectory (mixand) is selected from the working list $\mathcal{L}$, refined, and if successfully optimised, added to the solution set $\mathcal{L}^*$. The \gls*{DDP} solver is the same as the linear stochastic \gls*{AUL} solver of \Sec{sec:LSODA}.
\begin{algorithm}[h]
    \caption{Stochastic \gls*{DA}-based \gls*{AUL} solver}
    \label{alg:AUL_stochastic}
    \begin{algorithmic}
        \State \textbf{Input:} $\varepsilon_{\textrm{DDP}}$, $\varepsilon_{\textrm{AUL}}$, $\beta$, $\alpha_{\min}$, $\bm{x}_0$, $\bm{x}_t$, $\bm{U}_0$, $\bm{f}$, $l$, ${\phi}$, $\bm{g}$, $\bm{g}_t$, $\mathcal{T}$
        \State $\bm{U} \gets \bm{U}_0$
        \State Initialise $\bm{X}_0$
        \State $\beta^*, \ \mathcal{L}^*,\ \mathcal{L}\gets\beta,\ \emptyset,\ \left\{\left(\bm{X}_0,\bm{U}_0\right)\right\}$
        \While{$\mathcal{L}$ is not empty}
            \State Take the first pair $\left(\bm{X},\bm{U}\right)$ from $\mathcal{L}$
            \State Initialise $(\bm{X},\bm{U})$ using the element in $\mathcal{L}^*$ with the smallest Mahalanobis distance
            \State Initialise penalties and dual states
            \State Initialise augmented costs functions $\hat{l}$ and $\hat{\phi}$
            \State Compute $\mathcal{P}_{\bm{f}}$, $\mathcal{P}_{\mathcal{T}(\hat{l})}$, $\mathcal{P}_{\mathcal{T}(\hat{\phi})}$, $\bm{G}$, and $g_{\max}$
            \While{$g_{\max}>\varepsilon_{\textrm{AUL}}$}
                \State Compute $\eta\left(\bm{X},\bm{U},\mathcal{P}_{\bm{f}},\mathcal{P}_{\mathcal{T}(\hat{l})},\mathcal{P}_{\mathcal{T}(\hat{\phi})}\right)$
                \While{$\eta > \varepsilon_{\textrm{NLI}}$ and $\min\bm{\alpha}\geq\alpha_{\min}$}
                    \State $\left(\bm{X}_s,\bm{U}_s\right)\gets\text{GMMDecomposition}\left(\bm{X},\bm{U},\mathcal{P}_{\bm{f}},\mathcal{P}_{\mathcal{T}(\hat{l})},\mathcal{P}_{\mathcal{T}(\hat{\phi})}\right)$
                    \State Add $\left(\bm{X}_{-1},\bm{U}_{-1}\right)$ and $\left(\bm{X}_{1},\bm{U}_{1}\right)$ to the end of $\mathcal{L}$
                    \State $\bm{X},\ \bm{U} \gets \bm{X}_0,\ \bm{U}_0$ 
                    \State Compute $\eta\left(\bm{X},\bm{U},\mathcal{P}_{\bm{f}},\mathcal{P}_{\mathcal{T}(\hat{l})},\mathcal{P}_{\mathcal{T}(\hat{\phi})}\right)$
                \EndWhile
                \State $\Tilde{J}^*,\ \bm{X}^*,\ \bm{U}^* \gets \text{DDPSolver}\left(\bm{X}, \bm{U}, \bm{x}_t, \mathcal{P}_{\bm{f}}, \mathcal{T}(\hat{l},\beta^*), \mathcal{T}(\hat{\phi},\beta^*)\right)$
            \EndWhile
            \State Add $\left(\bm{X}^*, \bm{U}^*\right)$ to the end of $\mathcal{L}^*$
            \State $\beta_{\textrm{T}}\gets \text{RiskEstimation}(\bm{X}^*, \bm{U}^*,\mathcal{P}_{\mathcal{T}(\hat{l})},\mathcal{P}_{\mathcal{T}(\hat{\phi})})$
            \State $\beta^*\gets\text{FailureAllocationMethod}\left(\beta,\beta_{\textrm{T}},\bm{\alpha}\right)$
        \EndWhile
        \State \textbf{Return:} $\mathcal{L}^*$
    \end{algorithmic}
\end{algorithm}
Similarly, the projected Newton solver can be adapted to the non-linear case. \Alg{alg:pn_method_stochastic} presents the stochastic Newton solver.
\begin{algorithm}[h]
    \caption{Stochastic Newton solver}
    \label{alg:pn_method_stochastic}
    \begin{algorithmic}
        \State \textbf{Input:} $\beta$, $\varepsilon_{\textrm{N}}$, $\mathcal{L}$, $\bm{x}_t$, $\bm{f}$, $l$, ${\phi}$, $\bm{g}$, $\bm{g}_t$, $\mathcal{T}$
        \State $\beta^*, \ \mathcal{L}^* \gets\beta, \ \emptyset$
        \While{$\mathcal{L}$ is not empty}
            \State Take the first pair $\left(\bm{X},\bm{U}\right)$ from $\mathcal{L}$
            \State $\left(\bm{X}^*,\bm{U}^*\right)\gets\text{LinearStohasticNewtonMethod}\left(\bm{X},\bm{U},\bm{x}_t,\bm{f},\mathcal{T}(l,\beta^*),\mathcal{T}(\phi,\beta^*),\bm{g},\bm{g}_t\right)$ 
            \State Add $\left(\bm{X}^*, \bm{U}^*\right)$ to the end of $\mathcal{L}^*$
            \State $\beta_{\textrm{T}}\gets \text{RiskEstimation}(\bm{X}^*, \bm{U}^*,\mathcal{T}(l,\beta^*),\mathcal{T}(\phi,\beta^*))$
            \State $\beta^*\gets\text{FailureAllocationMethod}\left(\beta,\beta_{\textrm{T}},\bm{\alpha}\right)$
        \EndWhile
        \State \textbf{Return:} $\mathcal{L}^*$
    \end{algorithmic}
\end{algorithm}
Once the optimisation process yields a set of trajectories $\mathcal{L}^*$, the resulting control policies can be used to steer a spacecraft with initial condition $\hat{\bm{x}}_0$ drawn from the initial state distribution. To do so, the closest trajectory mixand of index $j$ is selected according to the squared Mahalanobis distance: $j = \text{argmin}_{i} \left(\hat{\bm{x}}_0 - \Bar{\bm{x}}_{0,i} \right)^{\textrm{T}} \bm{\Sigma}_{\bm{x},0,i}^{-1} \left(\hat{\bm{x}}_0 - \Bar{\bm{x}}_{0,i} \right)$
The corresponding feedback control law is then applied to propagate the sample trajectory: $\hat{\bm{u}}_{k} = \Bar{\bm{u}}_{k,j} + \bm{K}_{k,j} \left(\Bar{\bm{x}}_{k,j} - \hat{\bm{x}}_{k} \right)$, $\hat{\bm{x}}_{k+1} = \bm{f}\left(\hat{\bm{x}}_{k}, \hat{\bm{u}}_{k}\right) + \hat{\bm{w}}_{\bm{x},k+1}$,
where $\hat{\bm{w}}_{\bm{x},k+1}$ is a realisation of the Gaussian navigation noise $\bm{w}_{\bm{x},k+1}$.

%% file: sections/5_application.tex
\section{Applications}
\label{sec:application}

This section presents several applications of the proposed stochastic optimisation framework in the context of space mission design. We begin by validating the linear solver introduced in \Sec{sec:LSODA} using
the fuel-optimal Earth-Mars transfer from \citet{LantoineRussell_2012_AHDDPAfCOCPP2A}, followed by the assessment of the non-Gaussian solver from \Sec{sec:soda}. Finally, we demonstrate the application of our approach to trajectory design in the \gls*{CR3BP}.

All numerical results presented in this section were obtained using the \gls*{SODA} solver\footnote{Available at: \url{https://github.com/ThomasClb/SODA.git} [last accessed \lastdate].}, which was fully developed in C++. The only external dependency of the solver is the \gls*{DACE} library\footnote{Available at: \url{https://github.com/dacelib/dace.git} [last accessed \lastdate].}, developed by Dinamica Srl for the European Space Agency (ESA) \cite{RasottoEtAl_2016_DASTfNUPiSD}. All experiments and runtime analyses were carried out on a regular laptop equipped with an Intel\textsuperscript{\textregistered} Core\textsuperscript{TM} i7-11850H CPU running at 2.50\,GHz.

\subsection{Validation of the linear stochastic solver}
The \gls*{AUL} solver was validated on a stochastic fuel-optimal Earth–Mars transfer problem, adapted from \citet{LantoineRussell_2012_AHDDPAfCOCPP2A}. The dynamics of the low-thrust two-body problem, the fuel-optimal stage cost, and the terminal cost are the same as in \citet{CalebEtAl_2025_APBCSfFOLTTO}, and the stage constraints are given by $\bm{g}_{\textrm{ineq}}(\bm{x},\bm{u}) = \left[ \|\bm{u}\|^2 - u_{\textrm{max}}^2, m_{\textrm{dry}} - m \right]^{\textrm{T}}$, where $m$ denotes the spacecraft's total mass. In the remainder of this article, the initial spacecraft mass is \SI{1000}{\kilo\gram}, $m_{\textrm{dry}}=$\SI{500}{\kilo\gram} is the spacecraft's dry mass, $u_{\textrm{max}}=$\SI{0.5}{\newton} is the maximum thrust, and $\textrm{ISP}=$\SI{2000}{\second} is the thruster's specific impulse. Terminal constraints follow the definition in \Eq{eq:terminal_constraints}.  
The system has $N_x = 7$ states and $N_u = 3$ controls, with the state vector $\bm{x}$ expressed as $\left[x, y, z, \dot{x}, \dot{y}, \dot{z}, m\right]^{\textrm{T}}$. The trajectory is discretised into $N = 40$ stages over a fixed time of flight (ToF) of \SI{348.79}{\day}. The initial and target distributions are Gaussian with means and standard deviations reported in \Tab{tab:ic_earth_mars}. Each “Std.” row corresponds to the component-wise standard deviation of the associated distribution in normalised units. These uncertainty magnitudes are chosen such that the \gls*{NLI} remains below \num{5e-4}, ensuring the validity of the weakly non-linear assumption. The navigation noise covariances used throughout this work are set to $\bm{\mathcal{Q}}_{\bm{x},k} = \bm{\Sigma}_{\bm{x},0}/10^4$. Normalisation units are the same as in \citet{CalebEtAl_2025_APBCSfFOLTTO}.
\begin{table}[h]
    \centering
    \caption{Departure and target state distributions for the linear stochastic fuel-optimal Earth-Mars transfer.}
    \footnotesize
    \begin{tabular}{c c c c c c c c c} 
        \toprule
        \textbf{State} & \textbf{Quantity} & \textbf{$x$ [\SI{}{\kilo\meter}]} & \textbf{$y$ [\SI{}{\kilo\meter}]} & \textbf{$z$ [\SI{}{\kilo\meter}]} & \textbf{$\dot{x}$ [\SI{}{\kilo\meter}/\SI{}{\second}]} & \textbf{$\dot{y}$ [\SI{}{\kilo\meter}/\SI{}{\second}]} & \textbf{$\dot{z}$ [\SI{}{\kilo\meter}/\SI{}{\second}]} & \textbf{$m$ [\SI{}{\kilo\gram}]} \\
        \midrule
        \multirow{2}{*}{$\bm{x}_0$} & Mean & \num{-140699693} & \num{-51614428} & \num{980} & \num{9.774596} & \num{-28.07828} & \num{4.337725e-4} & \num{1000} \\
          & Std. & \num{e-6} & \num{e-6} & \num{e-6} & \num{5e-7} & \num{5e-7} & \num{5e-7} & 0 \\
        \midrule
        \multirow{2}{*}{$\bm{x}_t$}  &Mean  & \num{-172682023} & \num{176959469} & \num{7948912} & \num{-16.427384} & \num{-14.860506} & \num{9.21486e-2} & – \\
         & Std. & \num{e-4} & \num{e-4} & \num{e-4} & \num{e-5} & \num{e-5} & \num{e-5} & – \\
        \bottomrule
    \end{tabular}
    \label{tab:ic_earth_mars}
\end{table}
In the remainder of this article, all tolerance values are set to: $ \varepsilon_{\textrm{AUL}} = 10^{-6}$, $ \varepsilon_{\textrm{DDP}} = 10^{-4}$, $ \varepsilon_{\textrm{N}} = 10^{-10}$.
To ensure a balance between accuracy and computational efficiency, we enforce $\varepsilon_{\textrm{AUL}} = \varepsilon_{\textrm{DA}}$ throughout this article. This choice guarantees that only sufficiently accurate constraints are imposed, avoiding over-constraining the optimisation.

\begin{figure}[h]
    \centering
    \begin{subfigure}{.4\textwidth}
        \centering
        \includegraphics[width=1\linewidth]{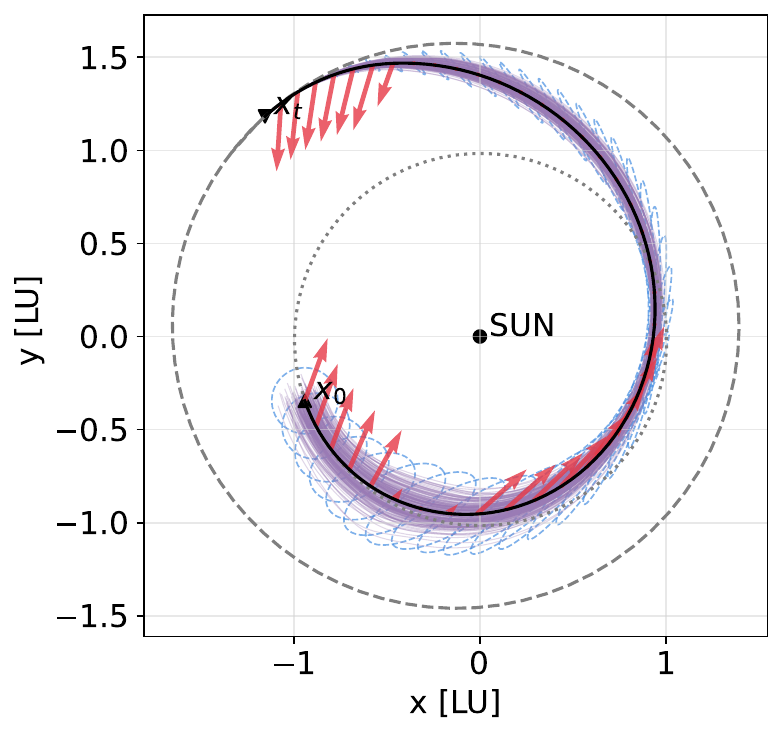}
        \caption{Trajectory in the $x\textrm{--}y$ plane.}
        \label{fig:earth_mars_L-SODA_trajectory}
    \end{subfigure}%
    \begin{subfigure}{.5\textwidth}
        \centering
        \includegraphics[width=1\linewidth]{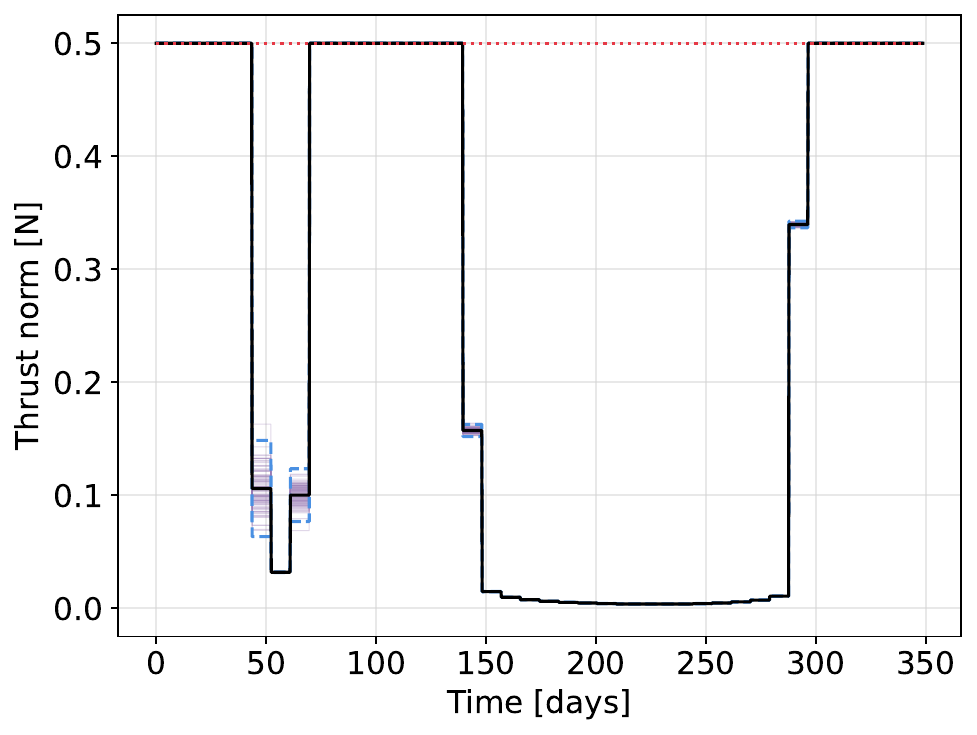}
        \caption{Thrust profile.}
        \label{fig:earth_mars_L-SODA_thrust}
    \end{subfigure}
    \caption{Solution to the linear stochastic fuel-optimal Earth-Mars transfer (errors magnified).}
    \label{fig:earth_mars_L-SODA}
\end{figure}

This stochastic optimisation problem was solved using two different approaches: \textbf{DET}: a deterministic formulation ignoring uncertainty. This solver corresponds to the \gls*{DADDy} solver from \citet{CalebEtAl_2025_APBCSfFOLTTO} and the proposed \gls*{L-SODA} solver using the first-order transcription method.
\Fig{fig:earth_mars_L-SODA} presents the solution resulting from L-SODA. The nominal Earth-to-Mars trajectory is shown in black in \Fig{fig:earth_mars_L-SODA_trajectory}, with red arrows representing the corresponding thrust vectors. In \Fig{fig:earth_mars_L-SODA_thrust}, the black line denotes the nominal control norm, the red dotted line indicates the maximum thrust $u_{\max}$, and the blue dashed curves represent theoretical \SI{95}{\%} confidence intervals. Several stochastic trajectory samples are also plotted in purple. All deviations were magnified by a factor $10^4$ for visualisation.
Quantitative results comparing the three approaches are shown in \Tab{tab:earth_mars_L-SODA_transcription}, for a target failure risk $\beta = 5\%$. All quantities ($1-\beta$ cost, failure rate $\beta_{\textrm{R}}$, and conservatism) were computed using Monte-Carlo sampling with \num{e5} realisations \citep{RobertCasella_2004_MCSM}.
\begin{table}[h]
    \centering
    \caption{Convergence data for the linear stochastic Earth--Mars transfer ($\beta = 5\%$).}
    \begin{tabular}{l c c c c c c c} 
        \toprule
        \textbf{Solver} & \textbf{RT [\SI{}{\second}]} & \textbf{\gls*{DDP} $n$}  & \textbf{Nominal error} & \textbf{$1-\beta$ cost [\SI{}{\kilo\gram}]} & \textbf{$\beta_{\textrm{R}}$ [\%]} & \textbf{$\gamma$} \\
        \midrule
        DET & \num{45.9} & \num{373} & \num{4e-13}  & \num{396.9} & \num{99.6} & \textendash \\
        L-SODA & \num{49.9} & \num{370} & \num{-1.22e-10} & \num{396.9} & \num{3.34} & \num{1.50} \\
        \bottomrule
    \end{tabular}
    \small
    \label{tab:earth_mars_L-SODA_transcription}
\end{table}
The deterministic solution (DET) fails to ensure constraint satisfaction under uncertainty, resulting in a \SI{99.6}{\%} failure rate. In contrast, the stochastic approach maintains feasibility for a similar number of iterations. In addition, L-SODA matches the deterministic solution in terms of runtime, iteration count, and $1{-}\beta$ cost, while achieving a low conservatism (\num{1.5}) and a failure rate of \SI{3.34}{\%}. This highlights the limited overhead induced by the proposed linear stochastic optimisation method (L-SODA) relative to deterministic solving (DET).
These results confirm the ability of the L-SODA solver to efficiently address weakly non-linear stochastic optimal control problems, with minimal added cost and conservatism. 

We now focus on the L-SODA solver and investigate its behaviour for various values of the target risk level $\beta$. The corresponding numerical results are presented in \Tab{tab:earth_mars_L-SODA_beta}. 
\begin{table}[h]
    \centering
    \caption{Convergence data for the L-SODA solver on the linear stochastic Earth--Mars transfer for various values of $\beta$.}
    \begin{tabular}{l c c c c c c c} 
        \toprule
        \textbf{$\beta$ [\%]} & \textbf{RT [\SI{}{\second}]} & \textbf{\gls*{DDP} $n$}  & \textbf{Nominal error} & \textbf{$1-\beta$ cost [\SI{}{\kilo\gram}]} & \textbf{$\beta_{\textrm{R}}$ [\%]} & \textbf{$\gamma$} \\
        \midrule
        \num{100} (DET) & \num{51.8} & \num{372} & \num{5e-15}  & \num{396.9} & \num{99.58} & \textendash \\
        \num{50} & \num{59.4} & \num{370} & \num{-1.79e-10} & \num{396.9} & \num{11.79} & \num{4.86} \\
        % \num{20} & \num{63.6} & \num{370} & \num{-9.27e-11} & \num{396.9} & \num{15.50} & \num{1.30} \\
        %\num{10} & \num{57.4} & \num{370} & \num{-1.09e-10} & \num{396.9} & \num{7.52} & \num{1.33} \\
        \num{5} & \num{49.9} & \num{370} & \num{-1.22e-10} & \num{396.9} & \num{3.34} & \num{1.50} \\
        %\num{2} & \num{61.6} & \num{370} & \num{-1.36e-10} & \num{396.9} & \num{1.11} & \num{1.81} \\
        %\num{1} & \num{61.4} & \num{370} & \num{-1.47e-10} & \num{396.9} & \num{0.50} & \num{1.98} \\
        \num{0.5} & \num{52.9} & \num{374} & \num{-1.32e-10} & \num{396.9} & \num{0.23} & \num{2.17} \\
        %\num{0.2} & \num{65.9} & \num{371} & \num{-1.41e-10} & \num{396.9} & \num{0.09} & \num{2.22} \\
        %\num{0.1} & \num{61.5} & \num{371} & \num{-1.47e-10} & \num{396.9} & \num{0.046} & \num{2.17} \\
        \bottomrule
    \end{tabular}
    \small
    \label{tab:earth_mars_L-SODA_beta}
\end{table}
Overall, the solver converges in a similar number of iterations across all considered values of $\beta$, and the resulting $1{-}\beta$ costs remain constant. Some variations in runtime are observed but remain within acceptable bounds and are likely due to numerical effects.
The real risk $\beta_{\textrm{R}}$ closely follows the target risk $\beta$, with conservatism values consistently below \num{5}, demonstrating robust and reliable behaviour of the method in all tested configurations.

\subsubsection{The stochastic fuel-optimal Earth-Mars transfer}
We now tackle the stochastic fuel-optimal Earth-Mars transfer without the weakly linear hypothesis, the new departure state covariance is $\bm{\Sigma}_{\bm{x},0}=\text{diag}\left[10^{-5},10^{-5},10^{-5},10^{-4},10^{-4},10^{-4},0\right]^2$. In the remainder of this paper, the target state distribution is $\bm{\Sigma}_t=\bm{\Sigma}_{\bm{x},0}/100$. We first study in the impact of the minimum splitting depth $\alpha_{\min}$ on the solutions. 
\Tab{tab:earth_mars_SODA_alpha} presents the convergence data of the SODA solver for several values of the minimum splitting depth $\alpha_{\min}$: \SI{50}{\%} no \gls*{GMM} is allowed as for the L-SODA solver, \SI{5}{\%}, and \SI{0.5}{\%}. The quantity $M$ stands for the number of mixands used in the \gls*{GMM}.
\begin{table}[h]
    \centering
    \caption{Convergence data for the SODA solver on the stochastic Earth--Mars transfer with $\beta=$\SI{5}{\%} for various values of $\alpha_{\min}$.}
    \begin{tabular}{l c c c c c c c c} 
        \toprule
        \textbf{$\alpha_{\min}$ [\%]} & \textbf{RT [\SI{}{\second}]} & \textbf{\gls*{DDP} $n$}  & \textbf{$M$} & \textbf{$\beta_{\textrm{T}}$ [\%]} & \textbf{$1-\beta$ cost [\SI{}{\kilo\gram}]}  & \textbf{$\beta_{\textrm{R}}$ [\%]} & \textbf{$\gamma$} \\
        \midrule
        \num{50} (L-SODA) & \num{269} & \num{315} & \num{1} & \num{0.36} &  \num{397.69} & \num{0.057} & \num{87.8} \\
        % \num{20} & \num{217.4} & \num{329} & \num{3} & \num{396.98} & \num{0.52} & \num{9.48} \\
        % \num{10} & \num{217.2} & \num{536} & \num{5} & \num{397.00} & \num{0.032} & \num{156} \\
        \num{5} & \num{273} & \num{574} & \num{11} & \num{3.13} &  \num{396.97} & \num{0.19} & \num{25.8} \\
        % \num{2} & \num{385.1} & \num{620} & \num{23} & \num{396.97} & \num{0.53} & \num{9.52} \\
        % \num{1} & \num{650.2} & \num{1606} & \num{45} & \num{396.97} & \num{0.21} & \num{24.2} \\
        \num{0.5} & \num{1030} & \num{3118} & \num{87} & \num{4.67} & \num{396.97} & \num{0.55} & \num{9.10} \\
        \bottomrule
    \end{tabular}
    \small
    \label{tab:earth_mars_SODA_alpha}
\end{table}
From these results, we observe that decreasing the minimum splitting depth $\alpha_{\min}$ leads, as expected, to an increase in the number of mixands. This, in turn, results in a higher number of \gls*{DDP} iterations and an overall increase in run time.  
However, the run time grows much more slowly than the number of mixands: for instance, while the number of mixands increases by a factor of \num{87}, the runtime is multiplied by less than \num{4}. This favourable scaling is primarily due to the educated initialisation strategy described in \Sec{sec:soda}.
In terms of feasibility, the estimated failure risk $\beta_{\textrm{T}}$ consistently satisfies $\beta_{\textrm{R}} \leq \beta_{\textrm{T}} \leq \beta$, and converges towards the target value $\beta$ more rapidly than the empirical risk $\beta_{\textrm{R}}$ as the number of mixands increases. From a performance standpoint, the SODA solver achieves a fuel saving of \SI{0.72}{\kilo\gram} compared to the L-SODA formulation.
Yet, \Tab{tab:earth_mars_SODA_alpha} shows no clear benefit in reducing $\alpha_{\min}$ below \SI{5}{\%}: while the $1{-}\beta$ cost remains essentially unchanged, the computational cost increases significantly.
As a result, we adopt $\alpha_{\min} = \SI{5}{\%}$ for the remainder of this work, as it offers a favourable trade-off between accuracy and computational efficiency. Although the computational cost of the SODA solver remains manageable across all tested settings, this configuration provides a significant improvement over L-SODA with limited additional overhead.

\begin{figure}[h]
    \centering
    \begin{subfigure}{.45\textwidth}
        \centering
        \includegraphics[width=1\linewidth]{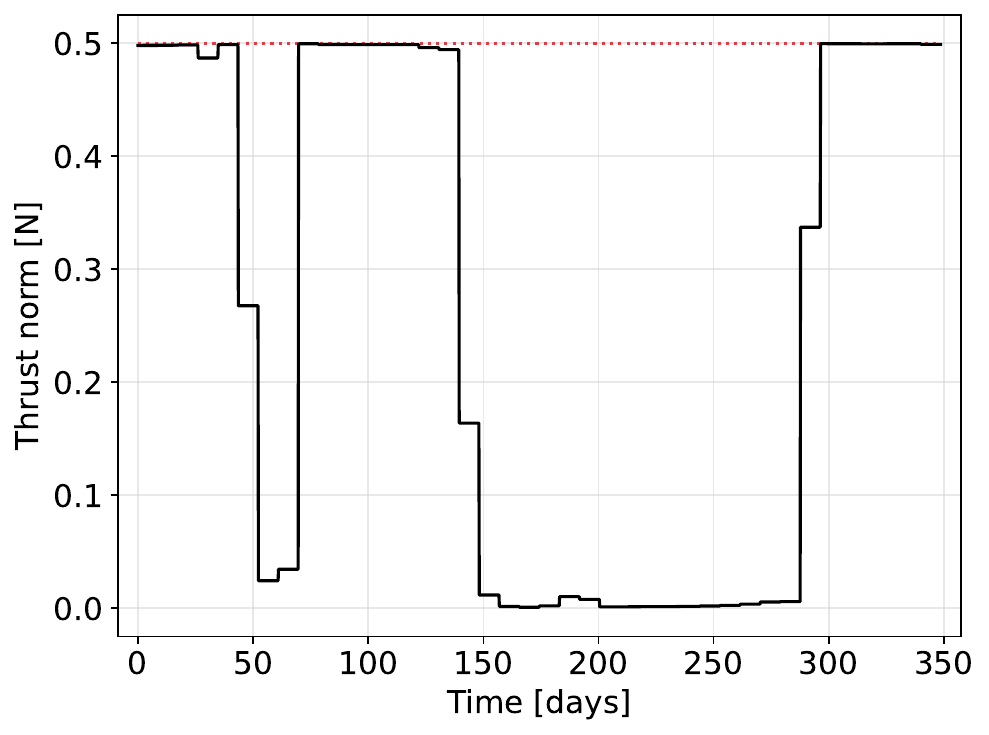}
        \caption{Thrust profile for $\alpha_{\min}=$\SI{50}{\%} (L-SODA).}
        \label{fig:earth_mars_SODA_trajectory_L-SODA}
    \end{subfigure}%
    \begin{subfigure}{.45\textwidth}
        \centering
        \includegraphics[width=1\linewidth]{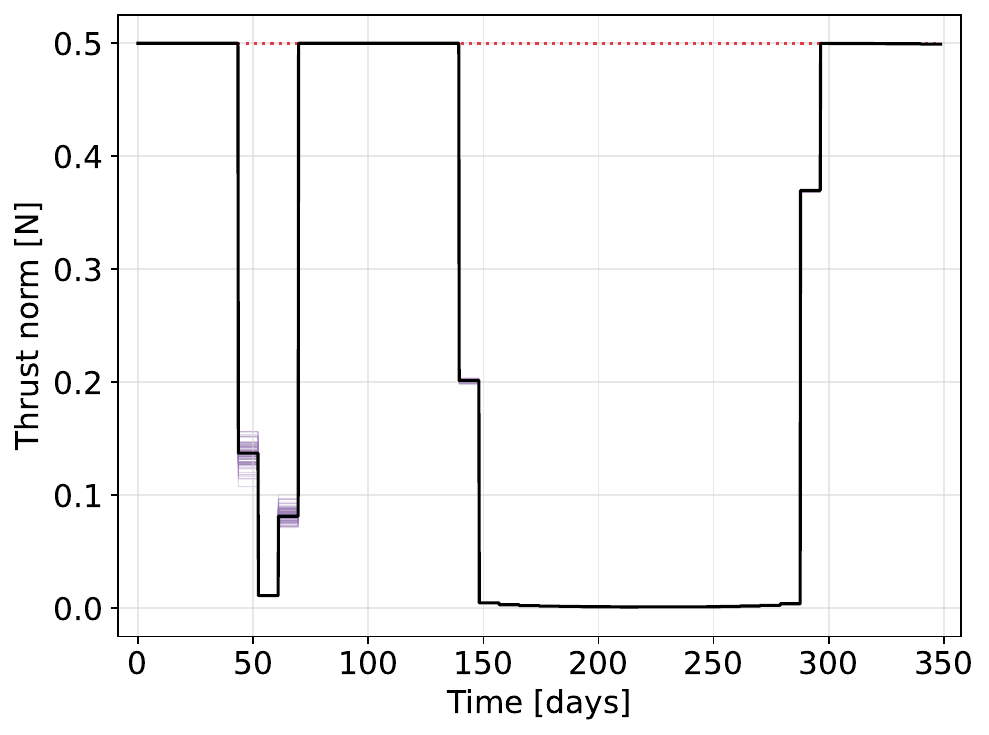}
        \caption{Thrust profile for $\alpha_{\min}=$\SI{5}{\%}.}
        \label{fig:earth_mars_SODA_thrust_SODA}
    \end{subfigure}
    \caption{Solutions to the stochastic fuel-optimal Earth-Mars transfer.}
    \label{fig:earth_mars_SODA_thrust}
\end{figure}

The gain in optimality offered by SODA over L-SODA can be illustrated in \Fig{fig:earth_mars_SODA_thrust}. \Fig{fig:earth_mars_SODA_trajectory_L-SODA} displays the thrust profile obtained with L-SODA, while \Fig{fig:earth_mars_SODA_thrust_SODA} shows the result for SODA with $\alpha_{\min} = \SI{5}{\%}$. In both figures, the black solid line denotes the nominal thrust, the red dotted line marks the thrust constraint, and the purple lines correspond to stochastic realisations of the thrust profile.
The L-SODA solution maintains large safety margins to account for uncertainty, which reduces optimality. In contrast, the SODA solver uses Gaussian mixture decomposition to optimise each mixand individually. This allows the solver to more tightly saturate the constraints, reducing fuel consumption.
Interestingly, the stochastic thrust realisations tend to cluster near the nominal profile when the nominal thrust is close to the maximum control magnitude. Conversely, greater variability appears in unconstrained phases, even without magnification. These deviations can be interpreted as correction manoeuvrers that exploit risk margins when the chance of violating the constraint is minimal.

\begin{figure}[h]
    \centering
    \begin{subfigure}{.5\textwidth}
        \centering
        \includegraphics[width=1\linewidth]{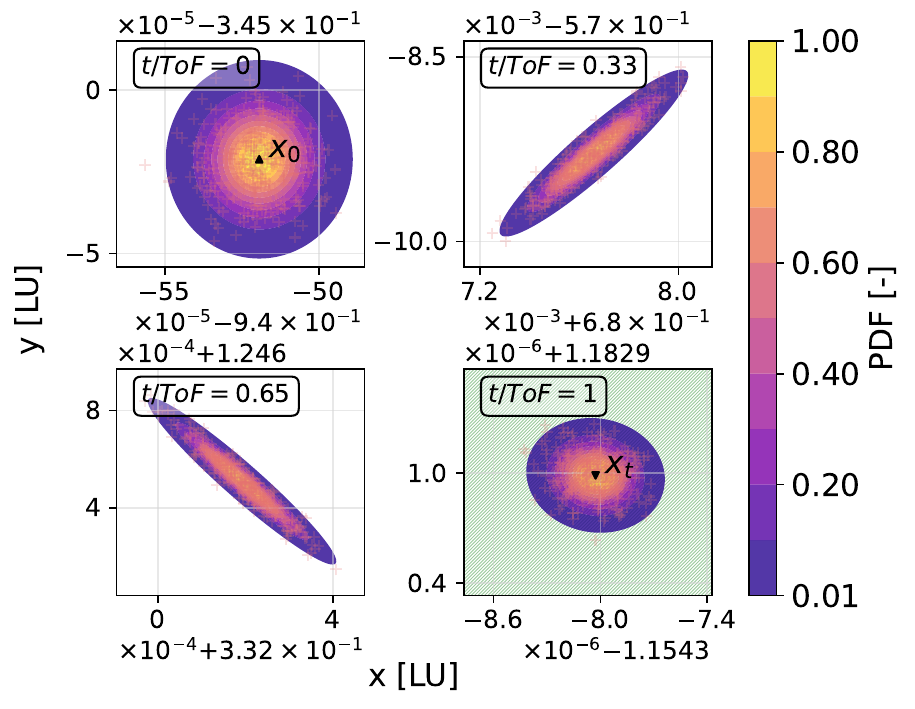}
        \caption{State \gls*{PDF} in the $x\textrm{--}y$ plane for $\alpha_{\min}=$\SI{50}{\%} (L-SODA).}
        \label{fig:earth_mars_SODA_pdf_position_L-SODA}
    \end{subfigure}%
    \begin{subfigure}{.5\textwidth}
        \centering
        \includegraphics[width=1\linewidth]{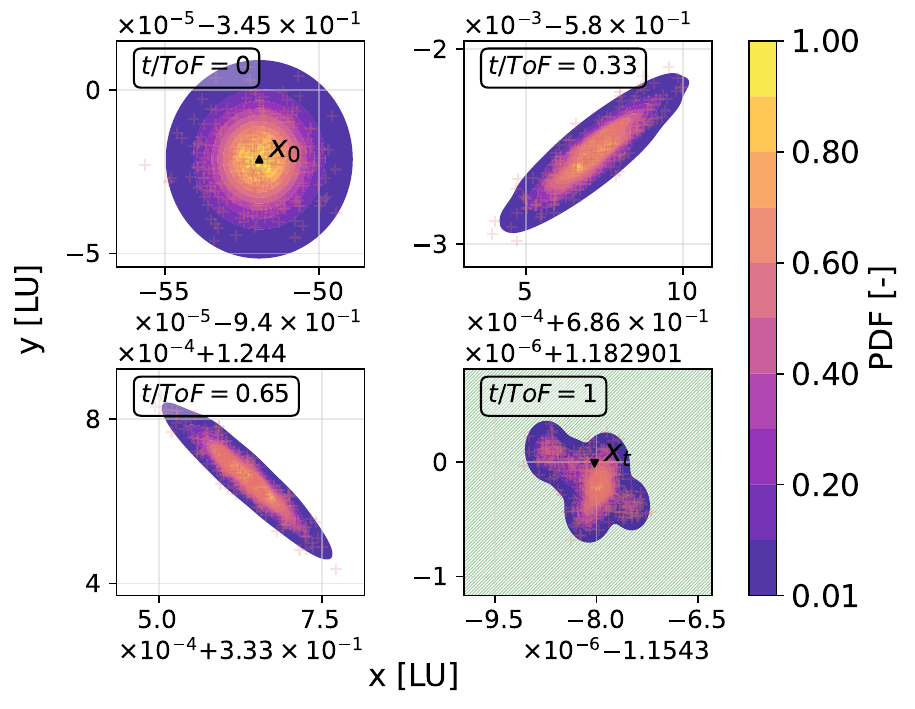}
        \caption{State \gls*{PDF} in the $x\textrm{--}y$ plane for $\alpha_{\min}=$\SI{5}{\%}.}
        \label{fig:earth_mars_SODA_pdf_position_SODA}
    \end{subfigure}
    \hfill
    \begin{subfigure}{.5\textwidth}
        \centering
        \includegraphics[width=1\linewidth]{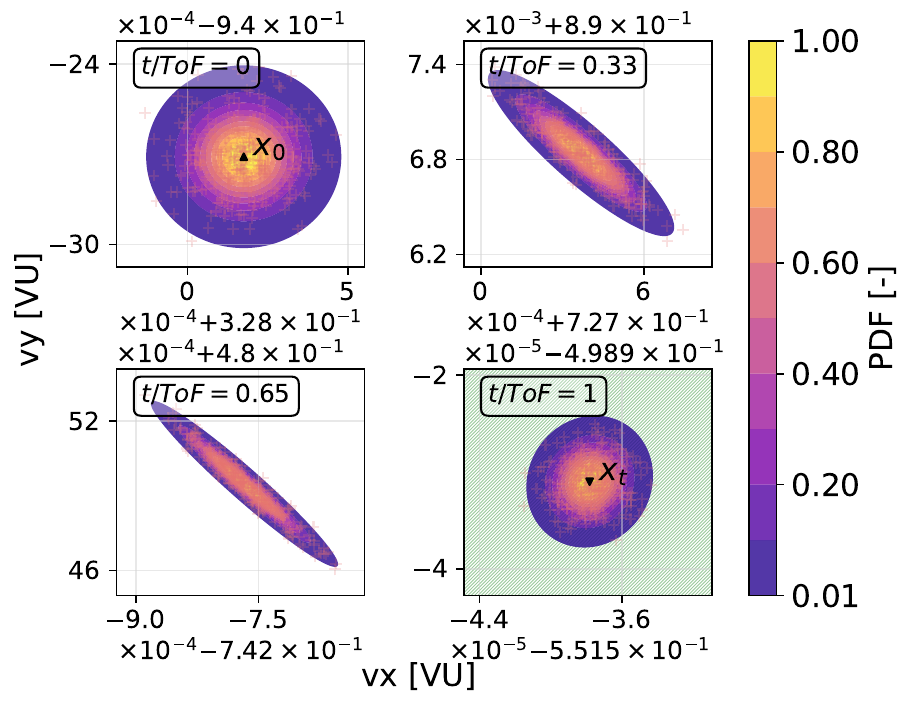}
        \caption{State \gls*{PDF} in the $\dot{x}\textrm{--}\dot{y}$ plane for $\alpha_{\min}=$\SI{50}{\%} (L-SODA).}
        \label{fig:earth_mars_SODA_pdf_velocity_L-SODA}
    \end{subfigure}%
    \begin{subfigure}{.5\textwidth}
        \centering
        \includegraphics[width=1\linewidth]{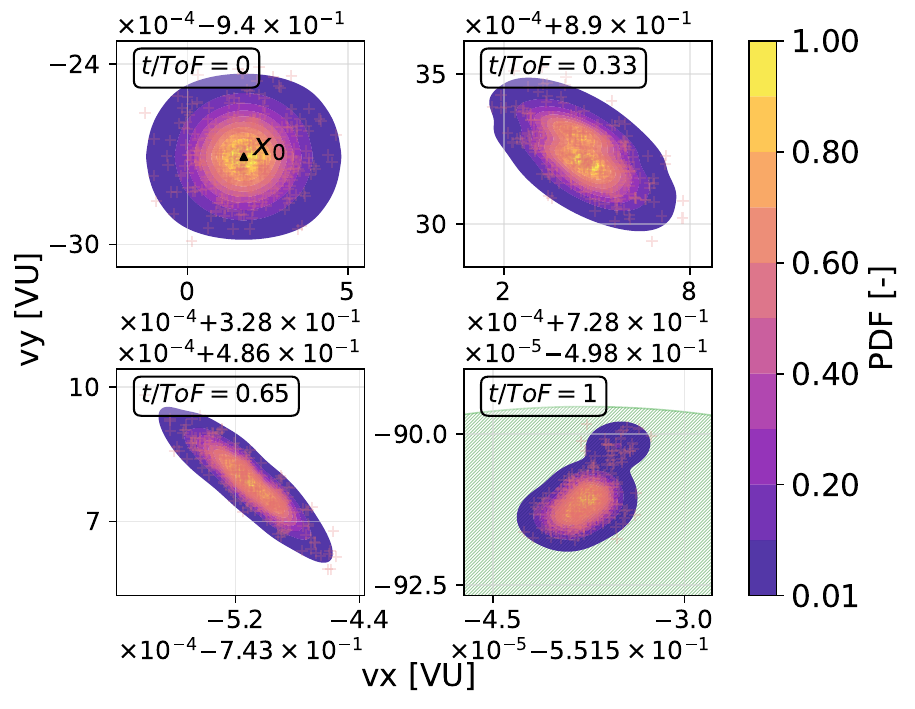}
        \caption{State \gls*{PDF} in the $\dot{x}\textrm{--}\dot{y}$ plane for $\alpha_{\min}=$\SI{5}{\%}.}
        \label{fig:earth_mars_SODA_pdf_velocity_SODA}
    \end{subfigure}
    \caption{Normalised state \gls*{PDF} of the stochastic fuel-optimal Earth-Mars transfer.}
    \label{fig:earth_mars_SODA_pdf}
\end{figure}

We now visualise the normalised \gls*{PDF}s of the position at different time steps in \Fig{fig:earth_mars_SODA_pdf}. \Fig{fig:earth_mars_SODA_pdf_position_L-SODA} and \Fig{fig:earth_mars_SODA_pdf_position_SODA} respectively show the normalised state \gls*{PDF} projected onto the $x\textrm{--}y$ plane for the L-SODA solver and the SODA solver with $\alpha_{\min} = \SI{5}{\%}$. Similarly, \Fig{fig:earth_mars_SODA_pdf_velocity_L-SODA} and \Fig{fig:earth_mars_SODA_pdf_velocity_SODA} display the corresponding distributions projected onto the $\dot{x}\textrm{--}\dot{y}$ plane. 
In each figure, the normalised theoretical \gls*{PDF} is shown along with stochastic realisations of the state, represented as orange crosses. The relative time corresponding to each window is indicated in the top-left corner; for instance, $t/\textrm{ToF} = 0$ corresponds to the initial state distribution, and $t/\textrm{ToF} = 1$ to the terminal state. In the latter case, the target area defined in \Eq{eq:terminal_constraints} is shown as a green hatched region.
These figures clearly indicate that both the L-SODA and SODA solvers are capable of steering the covariance with an accurate representation of the actual uncertainties, as the samples match the plotted \gls*{PDF}s. However, it is apparent that, in the case of the SODA solver, the mixands composing the distribution gradually separate over time, each following a distinct trajectory.

To complete the study of this test case, we now observe the behaviour of the \gls*{SODA} solver when the target failure risk varies in \Tab{tab:earth_mars_SODA_beta}.
\begin{table}[]
    \centering
    \caption{Convergence data for the SODA solver on the stochastic Earth--Mars transfer with $\alpha_{\min}=$\SI{5}{\%} for various values of $\beta$.}
    \begin{tabular}{l c c c c c c c c} 
        \toprule
        \textbf{$\beta$ [\%]} & \textbf{RT [\SI{}{\second}]} & \textbf{\gls*{DDP} $n$}  & \textbf{$M$} & \textbf{$\beta_{\textrm{T}}$ [\%]} & \textbf{$1-\beta$ cost [\SI{}{\kilo\gram}]} & \textbf{$\beta_{\textrm{R}}$ [\%]} & \textbf{$\gamma$} \\
        \midrule
        % \num{100} (DET) & \num{87.5} & \num{272} & \num{1} & \num{100} &  \num{396.96} & \num{99.85} & - \\
        \num{50} & \num{262} & \num{599} &  \num{9} & \num{36.86} &  \num{396.73} & \num{23.63} & \num{2.4} \\
        \num{5} & \num{273} & \num{574} & \num{11} & \num{3.13} &  \num{396.97} & \num{0.194} & \num{25.8} \\
        \num{0.5} & \num{386} & \num{582} & \num{11} & \num{0.12} & \num{397.68} & \num{0.048} & \num{10.4} \\
        \bottomrule
    \end{tabular}
    \small
    \label{tab:earth_mars_SODA_beta}
\end{table}
As expected, lowering the acceptable failure rate increases both computational effort and solution cost, since stricter chance constraints require larger safety margins. Notably, the conservatism exceeds \num{10} for $\beta = \SI{5}{\%}$ and $\beta = \SI{0.5}{\%}$, suggesting that further reductions in the $1{-}\beta$ cost may still be achievable. 
Indeed, the solution obtained for $\beta = \SI{5}{\%}$ already satisfies the constraints required for $\beta = \SI{0.5}{\%}$, indicating that the additional fuel consumption in the latter case is, \textit{a posteriori}, unnecessary. However, within the current formulation, there is no mechanism to relax the $\beta$ parameter \textit{a priori} to enable tighter constraint satisfaction without incurring additional cost. In particular, the estimated risk $\beta_{\textrm{T}}$ remains systematically closer to the target level $\beta$ than to the empirical failure rate $\beta_{\textrm{R}}$, thereby limiting the ability of the adaptive risk allocation strategy to relax the individual risk levels $\beta^*_j$ when appropriate.

\subsection{Test cases in the \gls*{CR3BP}}
We now tackle stochastic optimal control problems in the \gls*{CR3BP} \citep{Poincare_1892_LMNDLMC}. This dynamical system is famous for exhibiting strong non-linearities and provides important test cases to assess the performance of the \gls*{SODA} solver. The dynamics of the low-thrust \gls*{CR3BP} and the normalisation method are the same as in \citet{CalebEtAl_2025_APBCSfFOLTTO}, while the constraints, cost functions, spacecraft parameters, and tolerances are the same as in the Earth–Mars test case presented earlier.  

Similarly to \citet{CalebEtAl_2025_APBCSfFOLTTO}, the \gls*{SODA} solver was tested on three \gls*{CR3BP} test cases:
\begin{enumerate}
    \item A transfer from an $L_2$ halo orbit \cite{Farquhar_1970_TCaUoLPS} to an $L_1$ halo orbit used by \citet{AzizEtAl_2019_HDDPitCRTBP}.
    \item A transfer from an $L_2$ \gls*{NRHO} \citep{ZimovanSpreenEtAl_2020_NRHOaNHPDSOSaRP} to a \gls*{DRO} \cite{Henon_1969_NEotRPV}.
    % \item A Lyapunov $L_1$ \cite{Henon_1969_NEotRPV} to Lyapunov $L_2$ transfer from \citet{AzizEtAl_2019_HDDPitCRTBP}.
    \item A \gls*{DRO}-to-\gls*{DRO} transfer from \citet{AzizEtAl_2019_HDDPitCRTBP}.
\end{enumerate}
The initial conditions, targets, \gls*{ToF}s, and number of stages for each transfer are given in \Tab{tab:ic_cr3bp}.
\begin{table}[]
    \centering
    \caption{Earth-Moon \gls*{CR3BP} transfers data.}
    \small
    \begin{tabular}{l c c c c c c c c c} 
        \toprule
        \textbf{Transfer} & \textbf{\gls*{ToF}} [\SI{}{\day}] & $N$ & \textbf{Quantity} & \textbf{$x$ [LU]} & \textbf{$y$ [LU]} & \textbf{$z$ [LU]} & \textbf{$\dot{x}$ [VU]} & \textbf{$\dot{y}$ [VU]} & \textbf{$\dot{z}$ [VU]} \\ [0.5ex] 
        \midrule
        \multirow{3}{*}{\parbox{2cm}{Halo $L_2$ to halo $L_1$}} & \multirow{3}{*}{\num{20}} & \multirow{3}{*}{\num{110}} & $\bm{x}_0$ Mean & \num{1.16080} & \num{0} & \num{-0.12270} & \num{0} & \num{-0.20768} & \num{0} \\
        & & & $\bm{x}_0$ Std. & \num{e-6} & \num{e-6} & \num{e-6} & \num{e-5} & \num{e-5} & \num{e-5} \\
        & & & $\bm{x}_t$ Mean & \num{0.84871} & \num{0} & \num{0.17389} & \num{0} & \num{0.26350} & \num{0} \\
        \midrule[0.5pt]
        \multirow{3}{*}{\parbox{2cm}{\Gls*{NRHO} $L_2$ to \gls*{DRO}}} & \multirow{3}{*}{\num{21.2}} & \multirow{3}{*}{\num{150}} & $\bm{x}_0$ Mean & \num{1.02197} & \num{0} & \num{-0.18206} & \num{0} & \num{-0.10314} & \num{0} \\
        & & & $\bm{x}_0$ Std. & \num{2e-6} & \num{2e-6} & \num{2e-6} & \num{e-5} & \num{e-5} & \num{e-5} \\
        & & & $\bm{x}_t$ Mean & \num{0.98337} & \num{0.25921} & \num{0} & \num{0.35134} & \num{-0.00833} & \num{0} \\
        \midrule[0.5pt]
        %\multirow{3}{*}{\parbox{2cm}{Lyapunov $L_1$ to Lyapunov $L_2$}} & \multirow{3}{*}{\num{12}} & \multirow{3}{*}{\num{300}} & $\bm{x}_0$ Mean & \num{0.85599} & \num{0.12437} & \num{0} & \num{0.09485} & \num{0.04411} & \num{0} \\
        %& & & $\bm{x}_0$ Std. & \num{5e-6} & \num{5e-6} & \num{5e-6} & \num{5e-5} & \num{5e-5} & \num{5e-5} \\
        %& & & $\bm{x}_t$ Mean & \num{1.09598} & \num{0.11526} & \num{0} & \num{0.03747} & \num{0.12674} & \num{0}   \\
        %\midrule[0.5pt]
        \multirow{3}{*}{\parbox{2cm}{\gls*{DRO} to \gls*{DRO}}} & \multirow{3}{*}{\num{17.5}} & \multirow{3}{*}{\num{100}} & $\bm{x}_0$ Mean & \num{1.17136} & \num{0} & \num{0} & \num{0} & \num{-0.48946} & \num{0} \\
        & & & $\bm{x}_0$ Std. & \num{e-6} & \num{e-6} & \num{e-6} & \num{5e-5} & \num{5e-5} & \num{5e-5} \\
        & & & $\bm{x}_t$ Mean & \num{1.30184} & \num{0} & \num{0} & \num{0} & \num{-0.64218} & \num{0} \\
        \bottomrule
    \end{tabular}
    \label{tab:ic_cr3bp}
\end{table}
These test cases were solved using three different solvers: \textbf{DET}: the deterministic solver without uncertainty modelling; \textbf{\textbf{L-SODA}}: the proposed stochastic optimisation with differential algebra solver in the linear hypothesis ($\alpha_{\min}=50\%$); \textbf{\textbf{SODA}}: the proposed stochastic optimisation with differential algebra solver with $\alpha_{\min}=5\%$.

\begin{figure}[h]
    \centering 
    \begin{subfigure}{0.35\textwidth}
      \includegraphics[width=\linewidth]{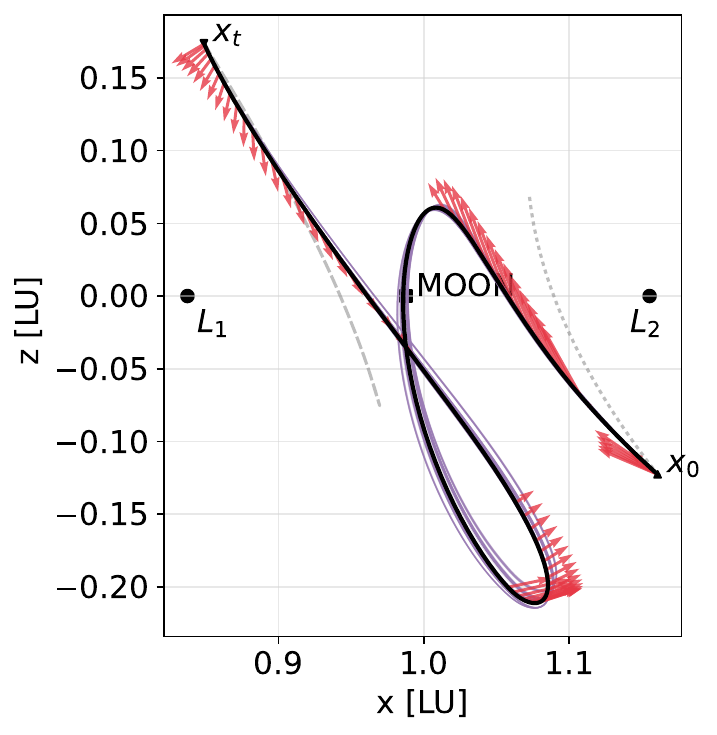}
      \caption{Trajectory in the $x\textbf{--}z$ plane (errors magnified).}
      \label{fig:halo_to_halo_xz}
    \end{subfigure}\hfil 
    \begin{subfigure}{0.47\textwidth}
      \includegraphics[width=\linewidth]{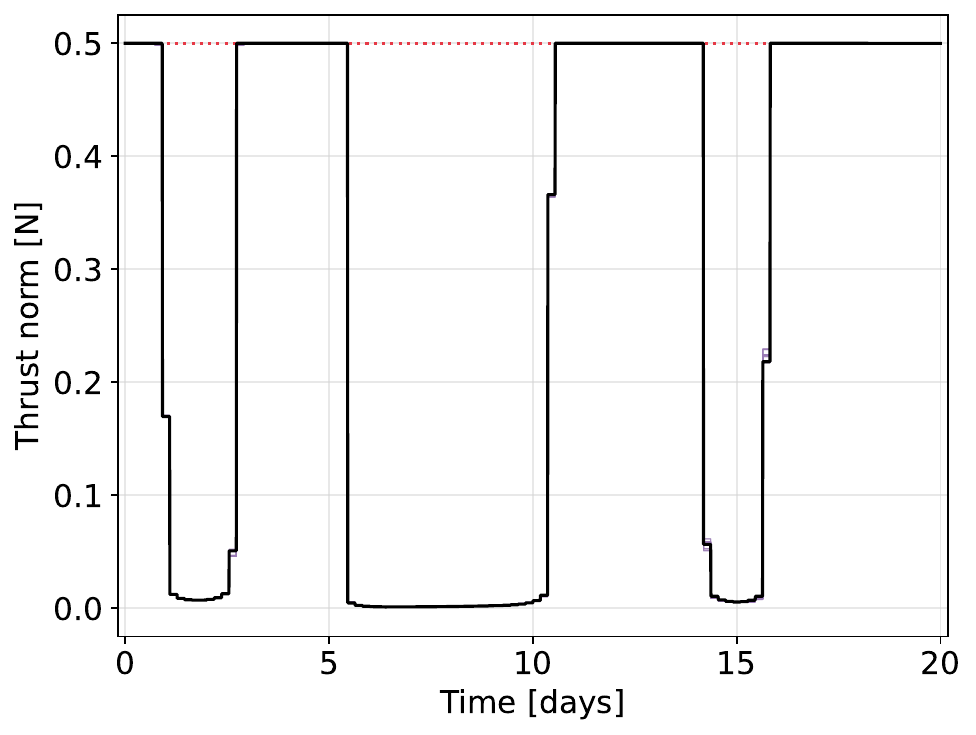}
      \caption{Thrust profile.}
      \label{fig:halo_to_halo_thrust}
    \end{subfigure}
    
    \medskip
    \begin{subfigure}{0.5\textwidth}
      \includegraphics[width=\linewidth]{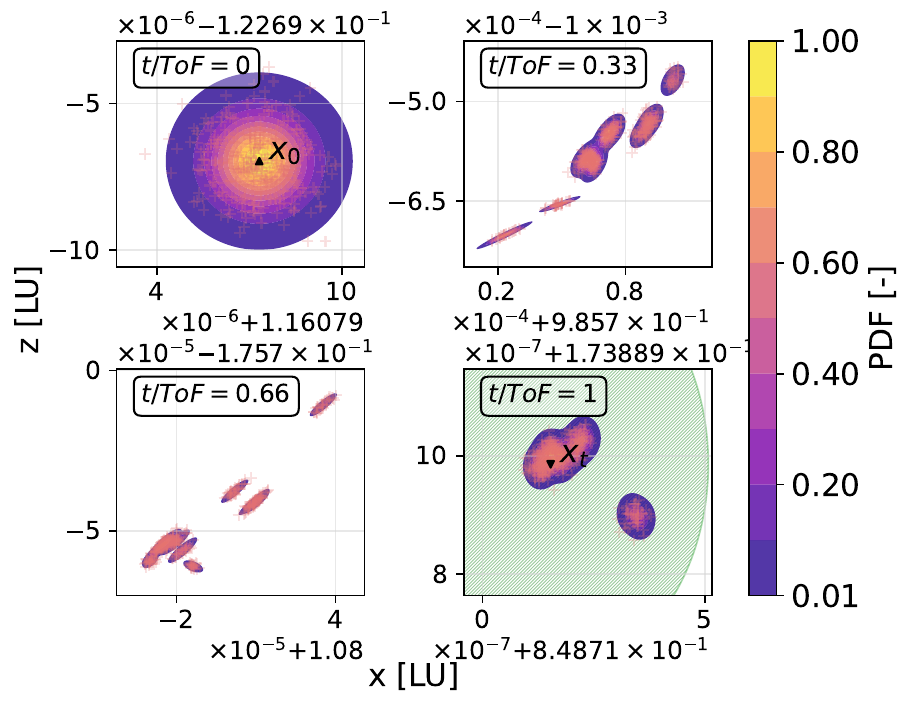}
      \caption{Normalised \gls*{PDF} in the $x\textbf{--}z$ plane.}
      \label{fig:halo_to_halo_pdf_position_xz}
    \end{subfigure}\hfil 
    \begin{subfigure}{0.5\textwidth}
      \includegraphics[width=\linewidth]{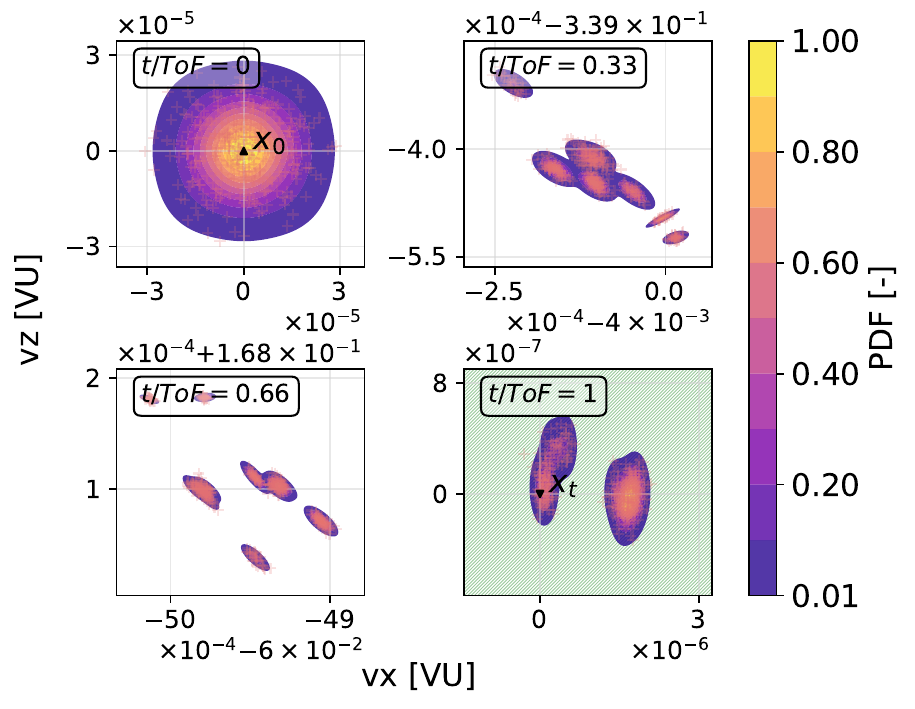}
      \caption{Normalised \gls*{PDF} in the $\dot{x}\textbf{--}\dot{z}$ plane.}
      \label{fig:halo_to_halo_pdf_velocity_xz}
    \end{subfigure}
    \caption{Solution to the stochastic halo $L_2$ to halo $L_1$ fuel-optimal transfer with SODA.}
    \label{fig:halo_to_halo}
\end{figure}

\begin{figure}[h]
    \centering
    \begin{subfigure}{.33\textwidth}
        \centering
        \includegraphics[width=1\linewidth]{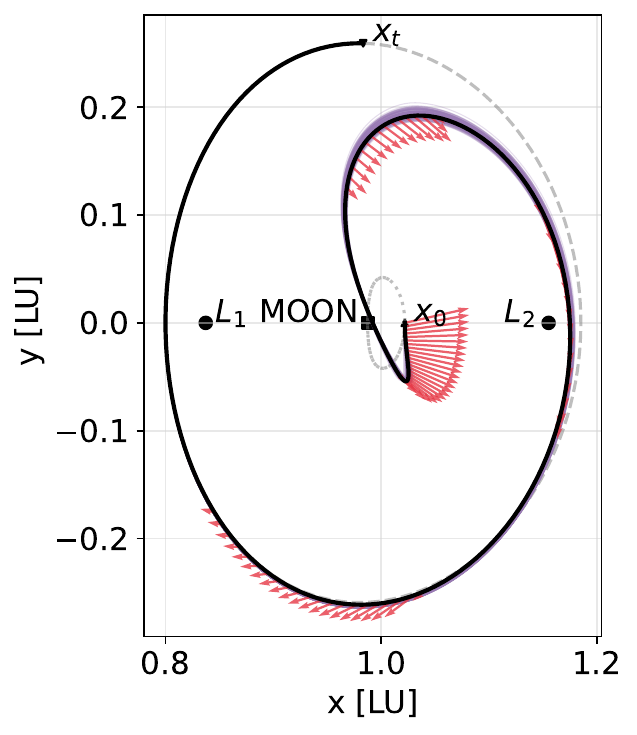}
        \caption{Trajectory in the $x\textbf{--}y$ plane (errors magnified).}
        \label{fig:nrho_to_dro_trajectory}
    \end{subfigure}%
    \begin{subfigure}{.51\textwidth}
        \centering
        \includegraphics[width=\linewidth]{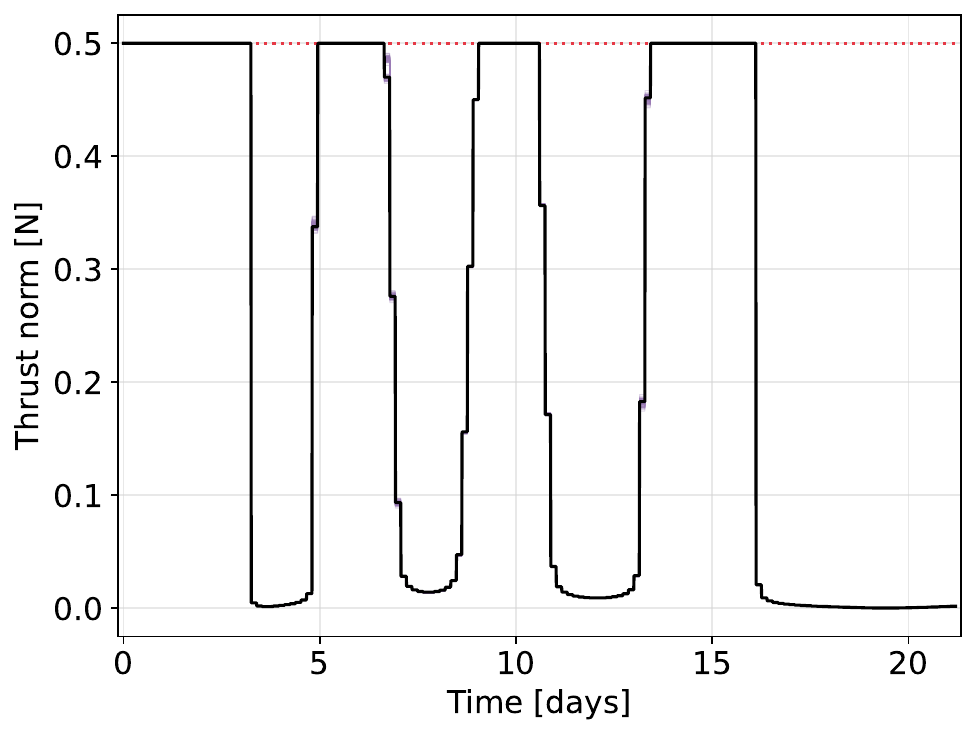}
        \caption{Thrust profile.}
        \label{fig:nrho_to_dro_thrust}
    \end{subfigure}
    \medskip
    \begin{subfigure}{0.5\textwidth}
      \includegraphics[width=\linewidth]{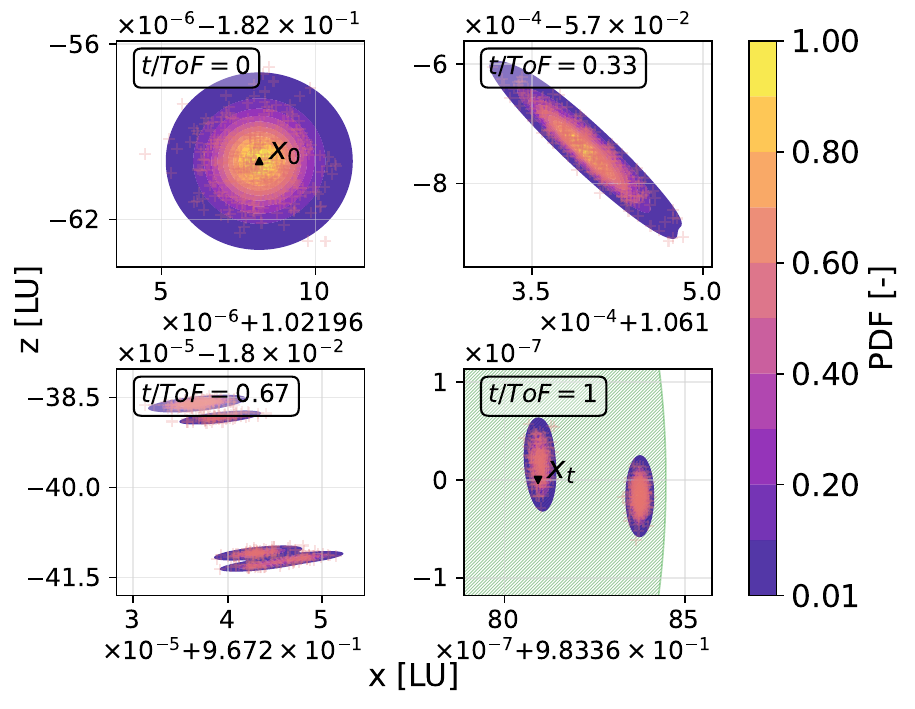}
      \caption{Normalised \gls*{PDF} in the $x\textbf{--}z$ plane.}
      \label{fig:nrho_to_dro_pdf_position_xz}
    \end{subfigure}\hfil 
    \begin{subfigure}{0.5\textwidth}
      \includegraphics[width=\linewidth]{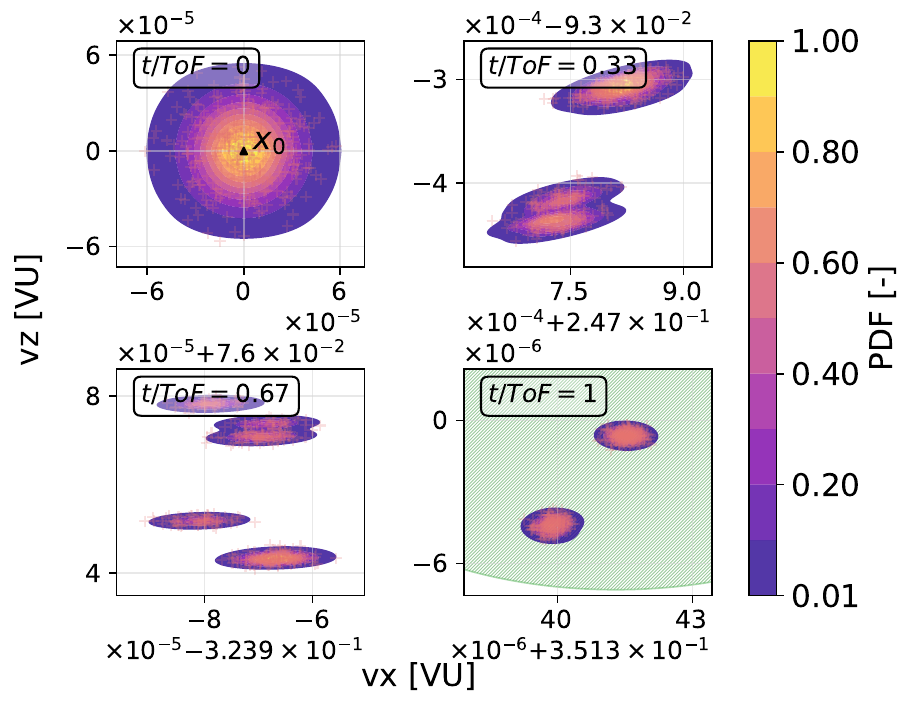}
      \caption{Normalised \gls*{PDF} in the $\dot{x}\textbf{--}\dot{z}$ plane.}
      \label{fig:nrho_to_dro_pdf_velocity_xz}
    \end{subfigure}
    \caption{Solution to the stochastic $L_2$ \gls*{NRHO} to \gls*{DRO} fuel-optimal transfer with SODA.}
    \label{fig:nrho_to_dro}
\end{figure}

\begin{figure}[h]
    \centering 
    \begin{subfigure}{0.29\textwidth}
      \includegraphics[width=0.85\linewidth]{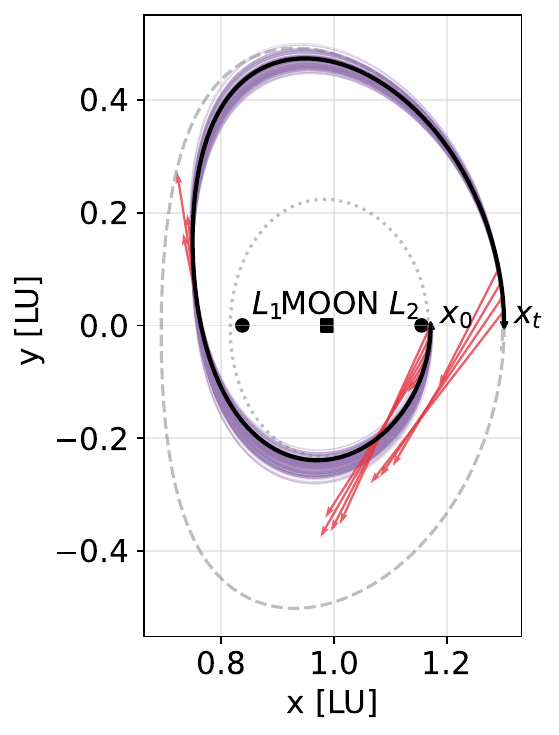}
      \caption{Trajectory in the $x\textbf{--}y$ plane (errors magnified).}
      \label{fig:dro_to_dro_trajectory}
    \end{subfigure}\hfil 
    \begin{subfigure}{0.55\textwidth}
      \includegraphics[width=0.85\linewidth]{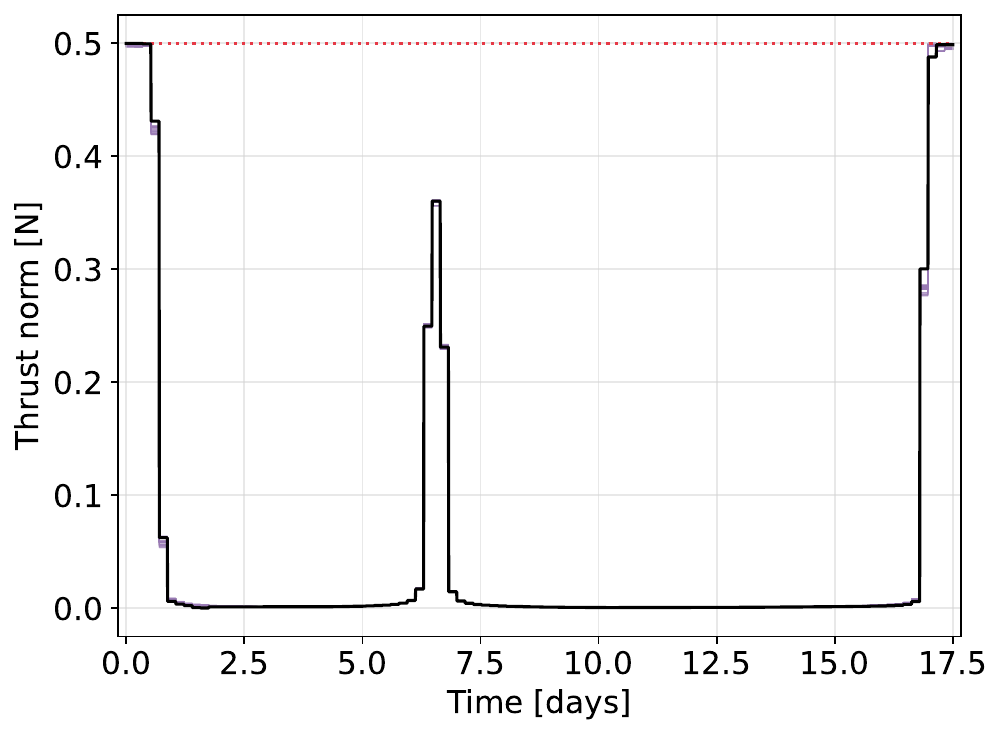}
      \caption{Thrust profile.}
      \label{fig:dro_to_dro_thrust}
    \end{subfigure}
    \medskip
    \begin{subfigure}{0.5\textwidth}
      \includegraphics[width=\linewidth]{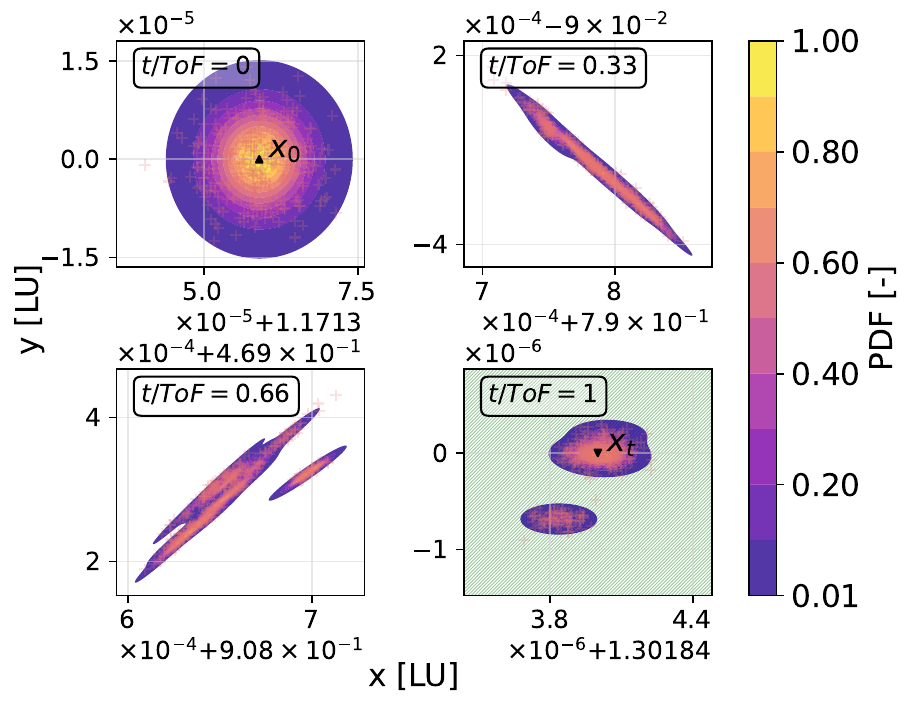}
      \caption{Normalised \gls*{PDF} in the $x\textbf{--}y$ plane.}
      \label{fig:dro_to_dro_pdf_position_xy}
    \end{subfigure}\hfil 
    \begin{subfigure}{0.5\textwidth}
      \includegraphics[width=\linewidth]{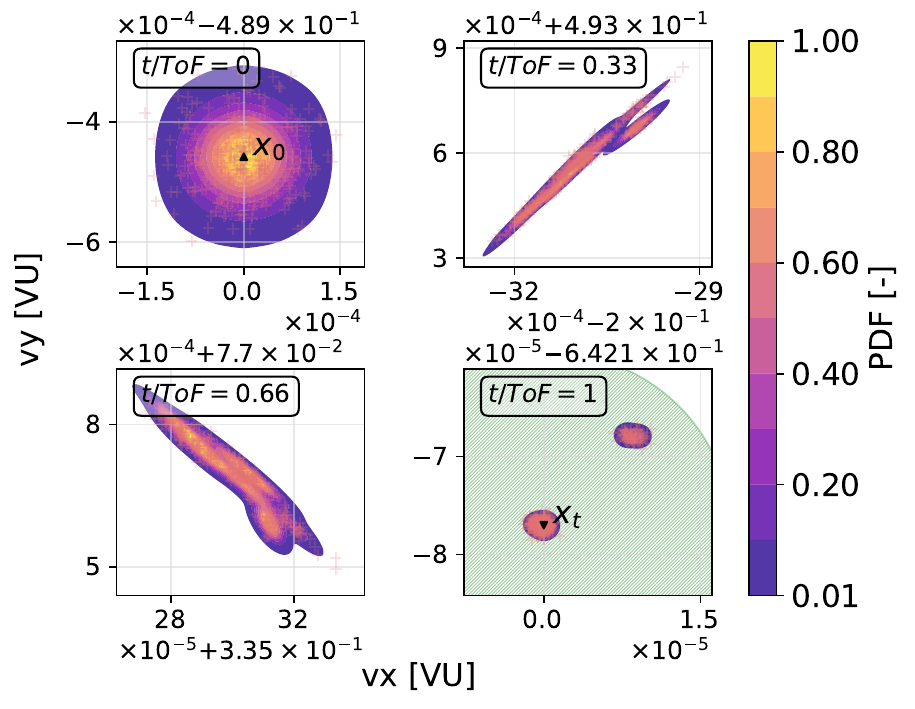}
      \caption{Normalised \gls*{PDF} in the $\dot{x}\textbf{--}\dot{y}$ plane.}
      \label{fig:dro_to_dro_pdf_velocity_xy}
    \end{subfigure}
    \caption{Solutions to the stochastic \gls*{DRO} to \gls*{DRO} fuel-optimal transfers with SODA.}
    \label{fig:dro_to_dro}
\end{figure}

\Fig{fig:halo_to_halo} shows the solutions to the stochastic halo $L_2$ to halo $L_1$ fuel-optimal transfers with SODA. \Fig{fig:halo_to_halo_xz}, \Fig{fig:halo_to_halo_thrust}, \Fig{fig:halo_to_halo_pdf_position_xz}, and \Fig{fig:halo_to_halo_pdf_velocity_xz} show respectively the transfer in the $x\textbf{--}z$ plane, the thrust profile, the normalised \gls*{PDF} in the $x\textbf{--}z$ plane, and the normalised \gls*{PDF} in the $\dot{x}\textbf{--}\dot{z}$ plane. 
\Fig{fig:nrho_to_dro} shows the solution to the stochastic \gls*{NRHO} to \gls*{DRO} fuel-optimal transfer with SODA. \Fig{fig:nrho_to_dro_trajectory}, \Fig{fig:nrho_to_dro_thrust}, \Fig{fig:nrho_to_dro_pdf_position_xz}, and \Fig{fig:nrho_to_dro_pdf_velocity_xz} show respectively the transfer in the $x\textbf{--}z$ plane, the thrust profile, the normalised \gls*{PDF} in the $x\textbf{--}z$ plane, and the normalised \gls*{PDF} in the $\dot{x}\textbf{--}\dot{z}$ plane.
Finally, \Fig{fig:dro_to_dro} shows the solution to the stochastic \gls*{DRO} to \gls*{DRO} fuel-optimal transfer with SODA. \Fig{fig:dro_to_dro_trajectory}, \Fig{fig:dro_to_dro_thrust}, \Fig{fig:dro_to_dro_pdf_position_xy}, and \Fig{fig:dro_to_dro_pdf_velocity_xy} show respectively the transfer in the $x\textbf{--}y$ plane, the thrust profile, the normalised \gls*{PDF} in the $x\textbf{--}y$ plane, and the normalised \gls*{PDF} in the $\dot{x}\textbf{--}\dot{y}$ plane.
These test cases exhibit qualitative behaviour similar to the Earth--Mars transfer: correction manoeuvrers are primarily located in regions where the thrust constraint is inactive, and the solver successfully contracts the distribution to satisfy the terminal constraints, even when the region of maximum covariance arises late in the trajectory, such as after the perilune. 
A key difference, however, lies in the earlier separation of the \gls*{PDF} into distinct trajectory bundles, leading to multiple well-separated branches rather than a single large trajectory fan. This behaviour, driven by stronger non-linearities, illustrates the solver’s ability to accommodate significant uncertainty growth while ensuring convergence toward the target set.

\Tab{tab:cr3bp_SODA} presents numerical results of the three solver on the \gls*{CR3BP} transfers.
\begin{table}[h]
    \centering
    \caption{Convergence data for the SODA solver on the stochastic \gls*{CR3BP} transfers with $\beta=$\SI{5}{\%}.}
    \small
    \begin{tabular}{l c c c c c c c c c} 
        \toprule
        \textbf{Transfer} & \textbf{Solver} & \textbf{RT [\SI{}{\second}]} & \textbf{\gls*{DDP} $n$}  & \textbf{$M$} & \textbf{$\beta_{\textrm{T}}$ [\%]} & \textbf{$1-\beta$ cost [\SI{}{\kilo\gram}]}  & \textbf{$\beta_{\textrm{R}}$ [\%]} & \textbf{$\gamma$} \\
        \midrule
        \multirow{3}{*}{\parbox{1.8cm}{Halo $L_2$ to halo $L_1$}} 
        & DET & \num{480} & \num{755} & \num{1} & \num{100} &  \num{26,068} & \num{99,98} & - \\
        & L-SODA & \num{494} & \num{347} & \num{1} & \num{0,45} &  \num{26,093} & \num{0,23} & \num{22} \\
        & SODA & \num{3004} & \num{2779} & \num{11} & \num{1,95} & \num{26,071} & \num{0,30} & \num{17} \\
        \midrule
        \multirow{3}{*}{\parbox{1.8cm}{\Gls*{NRHO} $L_2$ to \gls*{DRO}}} 
        & DET & \num{1061} & \num{531} & \num{1} & \num{100} &  \num{22,618} & \num{82,90} & - \\
        & L-SODA & \num{850} & \num{541} & \num{1} & \num{0} &  \num{22,650} & $<$\num{3.0e-3} & $>$\num{1700} \\
        & SODA & \num{3224} & \num{2010} & \num{11} & \num{2.61} & \num{22,619} & \num{2,44} & \num{2.1} \\
        \midrule
        %\multirow{3}{*}{\parbox{1.8cm}{Lyapunov $L_1$ to Lyapunov $L_2$}} 
        %& DET & \num{91} & \num{180} & \num{1} & \num{100} &  \num{3,703} & \num{98,53} & - \\
        %& L-SODA & \num{133} & \num{179} & \num{1} & \num{4,22} &  \num{3,726} & \num{0,63} & \num{8.0} \\
        %& SODA & \num{777} & \num{840} & \num{11} & \num{3,17} & \num{3,704} & \num{0,13} & \num{39} \\
        %\midrule
        \multirow{3}{*}{\parbox{2cm}{\gls*{DRO} to \gls*{DRO}}} 
        & DET & \num{122} & \num{267} & \num{1} & \num{100} &  \num{2,552} & \num{98,52} & - \\
        & L-SODA & \num{195} & \num{257} & \num{1} & \num{0} &  \num{2,624} & $<$\num{3.0e-3} & $>$\num{1700} \\
        & SODA & \num{1010} & \num{922} & \num{11} & \num{4,45} & \num{2,552} & \num{0,61} & \num{8.3} \\
        \bottomrule
    \end{tabular}
    \small
    \label{tab:cr3bp_SODA}
\end{table}
The deterministic (DET) and L-SODA solvers consistently exhibit similar run times and numbers of \gls*{DDP} iterations, reflecting their need to perform a single trajectory optimisation. In contrast, the SODA solver requires significantly more iterations due to the larger number of mixands (11). Nonetheless, the increase in computational cost remains moderate: the run time of SODA grows sub-linearly with respect to the number of mixands, and the ratio of run times between SODA and L-SODA remains well below this factor, highlighting the efficiency of the solver’s educated initialisation method.
From a feasibility standpoint, as anticipated, the DET solver fails to meet chance constraints, while both L-SODA and SODA consistently yield feasible solutions. This behaviour aligns with the theoretical guarantee that if $\beta_{\textrm{T}} \leq \beta$, then the trajectory is feasible and $\beta_{\textrm{R}} \leq \beta_{\textrm{T}}$. 
Note that cases with reported empirical failure rates below \SI{3.0e-3}{\percent} and conservatism above \num{1700} indicate no violations among \num{1e5} sampled trajectories, further supporting the important robustness of these solutions.
Moreover, except for the Lyapunov-to-Lyapunov transfer, the SODA solver is less conservative than L-SODA.
Regarding fuel consumption, the DET solver, as expected, achieves the lowest $1{-}\beta$ cost at the expense of feasibility. L-SODA enforces feasibility through conservative margins, resulting in noticeable sub-optimality. In contrast, the SODA solver maintains feasibility while significantly reducing the $1{-}\beta$ cost compared to L-SODA, with values approaching those of the DET reference. This demonstrates that SODA offers a compelling balance between robustness and performance.

%% file: sections/6_conclusion.tex
\section{Conclusion}
\label{sec:conclusion}
\glsresetall
In this paper, we presented \gls*{SODA}, a novel stochastic solver tailored to discrete-time chance-constrained space trajectory optimisation. Our method combines differential algebra with adaptive Gaussian mixture decomposition to efficiently propagate non-Gaussian uncertainties. Multidimensional Gaussian chance constraints are enforced using scalable and less conservative transcription methods, along with an efficient and accurate failure risk estimation strategy.

We validated the proposed framework through a series of test cases with increasing dynamical complexity, beginning with small uncertainty levels.
In a heliocentric transfer scenario, the \gls*{L-SODA} solver achieved performance comparable to a deterministic baseline, while providing robustness guarantees at negligible additional cost in computation and fuel.
Under more realistic uncertainty levels, \gls*{SODA} yielded improved constraint satisfaction and better optimality compared to \gls*{L-SODA}, even in strongly non-linear regimes such as the Earth–Moon \gls*{CR3BP}, where linearisation assumptions break down.
In these settings, \gls*{SODA} effectively handled non-Gaussian uncertainty propagation and multiple chance constrains, demonstrating a favourable balance between accuracy, robustness, and computational tractability.

Overall, \gls*{SODA} constitutes a coherent and effective framework for solving trajectory optimisation problems under uncertainty. Its ability to handle a broad class of dynamical systems and demonstrated performance across both linear and non-linear regimes suggest its strong potential for future applications in space mission analysis and design, particularly for uncertainty-aware trajectory planning in complex dynamical environments.

%% file: sections/7_funding_acknoledgments.tex
\section*{Funding sources}
This work was funded by SaCLaB (grant number 2022-CIF-R-1), a research group of ISAE-SUPAERO. 
\section*{Acknowledgments}
The authors acknowledge the use of an AI tool (Claude Haiku 4.5) to assist in the reformulation of certain sentences.
\section*{Declaration of competing interest}
The authors have no competing interests to declare that are relevant to the content of this article.

%% file: references.bib
@Article{LantoineRussell_2012_AHDDPAfCOCPP1T,
  author    = {Lantoine, Gregory and Russell, Ryan P.},
  journal   = {Journal of Optimization Theory and Applications},
  title     = {A Hybrid Differential Dynamic Programming Algorithm for Constrained Optimal Control Problems. Part 1: Theory},
  year      = {2012},
  month     = {4},
  number    = {2},
  pages     = {382--417},
  volume    = {154},
  doi       = {10.1007/s10957-012-0039-0},
  publisher = {Springer Science and Business Media {LLC}},
  url       = {https://doi.org/10.1007/s10957-012-0039-0},
}

@InProceedings{CaoEtAl_2002_PCFoaBTM,
  author    = {Cao, T.D. and Hall, J.F. and van de Geijn, R.A.},
  booktitle = {Proceedings. International Conference on Parallel Processing Workshop},
  title     = {Parallel Cholesky factorization of a block tridiagonal matrix},
  year      = {2002},
  pages     = {327-335},
  doi       = {10.1109/ICPPW.2002.1039748},
}

@InProceedings{HowellEtAl_2019_AaFSfCTO,
  author    = {Howell, Taylor A. and Jackson, Brian E. and Manchester, Zachary},
  booktitle = {2019 IEEE/RSJ International Conference on Intelligent Robots and Systems (IROS)},
  title     = {ALTRO: A Fast Solver for Constrained Trajectory Optimization},
  year      = {2019},
  pages     = {1-6},
  doi       = {10.1109/IROS40897.2019.8967788},
  url       = {https://doi.org/10.1109/IROS40897.2019.8967788},
}

@Article{BenedikterEtAl_2022_CAtCCwAtSLTTO,
  author  = {Benedikter, Boris and Zavoli, Alessandro and Wang, Zhenbo and Pizzurro, Simone and Cavallini, Enrico},
  journal = {Journal of Guidance, Control, and Dynamics},
  title   = {Convex Approach to Covariance Control with Application to Stochastic Low-Thrust Trajectory Optimization},
  year    = {2022},
  number  = {11},
  pages   = {2061-2075},
  volume  = {45},
  doi     = {10.2514/1.G006806},
  url     = {https://doi.org/10.2514/1.G006806},
}

@Article{LantoineRussell_2012_AHDDPAfCOCPP2A,
  author    = {Lantoine, Gregory and Russell, Ryan P.},
  journal   = {Journal of Optimization Theory and Applications},
  title     = {A Hybrid Differential Dynamic Programming Algorithm for Constrained Optimal Control Problems. Part 2: Application},
  year      = {2012},
  month     = {8},
  number    = {2},
  pages     = {418--442},
  volume    = {154},
  doi       = {10.1007/s10957-012-0038-1},
  publisher = {Springer Science and Business Media {LLC}},
  url       = {https://doi.org/10.1007/s10957-012-0038-1},
}

@Article{ArmellinEtAl_2010_ACECUDAtCoA,
  author    = {Armellin, R and Di Lizia, P and Bernelli-Zazzera, F and Berz, M},
  journal   = {Celestial Mechanics and Dynamical Astronomy},
  title     = {Asteroid close encounters characterization using differential algebra: the case of {Apophis}},
  year      = {2010},
  month     = {6},
  number    = {4},
  pages     = {451--470},
  volume    = {107},
  doi       = {10.1007/s10569-010-9283-5},
  publisher = {Springer Science and Business Media {LLC}},
}

@Article{WittigEtAl_2015_PoLUSiODbADS,
  author    = {Wittig, A and Di Lizia, P and Armellin, R and Makino, K and Bernelli-Zazzera, F and Berz, M},
  journal   = {Celestial Mechanics and Dynamical Astronomy},
  title     = {Propagation of large uncertainty sets in orbital dynamics by automatic domain splitting},
  year      = {2015},
  month     = {5},
  number    = {3},
  pages     = {239--261},
  volume    = {122},
  doi       = {10.1007/s10569-015-9618-3},
  publisher = {Springer Science and Business Media {LLC}},
}

@Article{CalebEtAl_2023_DAMAtCAGaBDAtPFotES,
  author    = {Caleb, T and Losacco, M and Fossà, A and Armellin, R and Lizy-Destrez, S},
  journal   = {Nonlinear Dynamics},
  title     = {{Differential Algebra Methods Applied to Continuous Abacus Generation and Bifurcation Detection: Application to Periodic Families of the Earth–Moon system}},
  year      = {2023},
  month     = {3},
  doi       = {10.1007/s11071-023-08375-0},
  publisher = {Springer Science and Business Media {LLC}},
}

@Article{Henon_1969_NEotRPV,
  author  = {Hénon, M},
  journal = {Astronomy and Astrophysics},
  title   = {Numerical exploration of the restricted problem, V},
  year    = {1969},
  pages   = {223-238},
  volume  = {1},
}

@Article{ZimovanSpreenEtAl_2020_NRHOaNHPDSOSaRP,
  author    = {Zimovan Spreen, E M and Howell, K C and Davis, D C},
  journal   = {Celestial Mechanics and Dynamical Astronomy},
  title     = {Near rectilinear halo orbits and nearby higher-period dynamical structures: orbital stability and resonance properties},
  year      = {2020},
  month     = {5},
  number    = {5},
  pages     = {1--25},
  volume    = {132},
  doi       = {10.1007/s10569-020-09968-2},
  publisher = {Springer Science and Business Media {LLC}},
}

@Article{AzizEtAl_2019_HDDPitCRTBP,
  author  = {Aziz, Jonathan D. and Scheeres, Daniel J. and Lantoine, Gregory},
  journal = {Journal of Guidance, Control, and Dynamics},
  title   = {Hybrid Differential Dynamic Programming in the Circular Restricted Three-Body Problem},
  year    = {2019},
  number  = {5},
  pages   = {963-975},
  volume  = {42},
  doi     = {10.2514/1.G003617},
  url     = {https://doi.org/10.2514/1.G003617},
}

@Article{Mayne_1966_ASOGMfDOToNLDTS,
  author    = {Mayne, David},
  journal   = {International Journal of Control},
  title     = {A Second-order Gradient Method for Determining Optimal Trajectories of Non-linear Discrete-time Systems},
  year      = {1966},
  number    = {1},
  pages     = {85--95},
  volume    = {3},
  doi       = {10.1080/00207176608921369},
  publisher = {Taylor \& Francis},
  url       = {https://doi.org/10.1080/00207176608921369},
}

@InProceedings{RasottoEtAl_2016_DASTfNUPiSD,
  author    = {Rasotto, M and Morselli, A and Wittig, A and Massari, M and Di Lizia, P and Armellin, R and Yabar Valles, C and Urbina Ortega, C},
  booktitle = {6th International Conference on Astrodynamics Tools and Techniques (ICATT)},
  title     = {Differential algebra space toolbox for nonlinear uncertainty propagation in space dynamics},
  year      = {2016},
  editor    = {uropean Space Operations Centre (ESOC), Darmstadt, Germany},
  month     = {3},
  pages     = {1--11},
}

@InProceedings{SmithEtAl_2020_TAPaOoNAtRHttM,
  author    = {Smith, M and Craig, D and Herrmann, N and Mahoney, E and Krezel, J and McIntyre, N and Goodliff, K},
  booktitle = {2020 IEEE Aerospace Conference},
  title     = {{The Artemis Program: An Overview of NASA's Activities to Return Humans to the Moon}},
  year      = {2020},
  pages     = {1-10},
  doi       = {10.1109/AERO47225.2020.9172323},
}

@InBook{Berz_1999_MMMiPBP,
  author    = {Berz, M},
  chapter   = {2},
  pages     = {82--119},
  publisher = {Elsevier},
  title     = {Modern Map Methods in Particle Beam Physics},
  year      = {1999},
  doi       = {10.1016/s1076-5670(08)70227-1},
}

@Book{Poincare_1892_LMNDLMC,
  author    = {Poincaré, H},
  publisher = {Gauthier-Villars, Paris, France},
  title     = {Les Méthodes Nouvelles de la Mécanique Céleste},
  year      = {1892},
  language  = {fra},
  pages     = {1--12},
}

@InBook{RobertCasella_2004_MCSM,
  author    = {Robert, C and Casella, G},
  chapter   = {3},
  pages     = {79--122},
  publisher = {Springer New York},
  title     = {Monte Carlo Statistical Methods},
  year      = {2004},
  doi       = {10.1007/978-1-4757-4145-2},
}

@TechReport{Farquhar_1970_TCaUoLPS,
  author      = {Farquhar, R W},
  institution = {National Aeronautics and Space Administration (NASA)},
  title       = {The Control and Use of Libration-Point Satellites},
  year        = {1970},
  note        = {number NASA TR R-346},
  pages       = {1--21},
  school      = {National Aeronautics and Space Administration (NASA)},
}

@TechReport{Cholesky_1910_SLRNDSDL,
  author = {Cholesky, André-Louis},
  title  = {Sur la résolution numérique des systèmes d'équations linéaires},
  year   = {1910},
}

@Article{OzakiEtAl_2018_SDDPwUTfLTTD,
  author  = {Ozaki, Naoya and Campagnola, Stefano and Funase, Ryu and Yam, Chit Hong},
  journal = {Journal of Guidance, Control, and Dynamics},
  title   = {Stochastic Differential Dynamic Programming with Unscented Transform for Low-Thrust Trajectory Design},
  year    = {2018},
  number  = {2},
  pages   = {377-387},
  volume  = {41},
  doi     = {10.2514/1.G002367},
  url     = {https://doi.org/10.2514/1.G002367},
}

@Article{OzakiEtAl_2020_TSOCfNCTOP,
  author  = {Ozaki, Naoya and Campagnola, Stefano and Funase, Ryu},
  journal = {Journal of Guidance, Control, and Dynamics},
  title   = {Tube Stochastic Optimal Control for Nonlinear Constrained Trajectory Optimization Problems},
  year    = {2020},
  number  = {4},
  pages   = {645-655},
  volume  = {43},
  doi     = {10.2514/1.G004363},
  url     = {https://doi.org/10.2514/1.G004363},
}

@InProceedings{BooneMcMahon_2022_NGCCTCUGMaRA,
  author    = {Boone, S. and McMahon, J.},
  booktitle = {IEEE 61st Conference on Decision and Control (CDC)},
  title     = {Non-Gaussian Chance-Constrained Trajectory Control Using Gaussian Mixtures and Risk Allocation},
  year      = {2022},
  pages     = {3592--3597},
  ranking   = {rank2},
}

@InProceedings{RidderhofEtAl_2020_CCCCfLTMFTO,
  author    = {Ridderhof, J. and Pilipovsky, J. and Tsiotras, P.},
  booktitle = {2020 {AIAA}/{AAS} Astrodynamics Specialist Conference},
  title     = {Chance-Constraints Covariance Control for Low-Thrust Minimum-Fuel Trajectory Optimization},
  year      = {2020},
  pages     = {1--20},
}

@Article{Mahalanobis_1936_OtGDiS,
  author  = {Mahalanobis, P. C.},
  journal = {Proceedings of the National Institute of Sciences of India},
  title   = {On the Generalized Distance in Statistics},
  year    = {1936},
  number  = {1},
  pages   = {49-55},
  volume  = {2},
  url     = {https://www.jstor.org/stable/10.2307/48723335},
}

@Article{LosaccoEtAl_2024_LOADSAfNUM,
  author    = {Losacco, Matteo and Fossà, Alberto and Armellin, Roberto},
  journal   = {Journal of Guidance, Control, and Dynamics},
  title     = {Low-Order Automatic Domain Splitting Approach for Nonlinear Uncertainty Mapping},
  year      = {2024},
  issn      = {1533-3884},
  month     = feb,
  number    = {2},
  pages     = {291--310},
  volume    = {47},
  doi       = {10.2514/1.g007271},
  publisher = {American Institute of Aeronautics and Astronautics (AIAA)},
}

@Article{DeMarsEtAl_2013_EBAfUPoNDS,
  author    = {DeMars, Kyle J. and Bishop, Robert H. and Jah, Moriba K.},
  journal   = {Journal of Guidance, Control, and Dynamics},
  title     = {Entropy-Based Approach for Uncertainty Propagation of Nonlinear Dynamical Systems},
  year      = {2013},
  issn      = {1533-3884},
  month     = jul,
  number    = {4},
  pages     = {1047--1057},
  volume    = {36},
  doi       = {10.2514/1.58987},
  publisher = {American Institute of Aeronautics and Astronautics (AIAA)},
}

@Article{BlackmoreEtAl_2010_APPCAoCCSPC,
  author    = {Blackmore, Lars and Ono, Masahiro and Bektassov, Askar and Williams, Brian C.},
  journal   = {IEEE Transactions on Robotics},
  title     = {A Probabilistic Particle-Control Approximation of Chance-Constrained Stochastic Predictive Control},
  year      = {2010},
  issn      = {1941-0468},
  month     = jun,
  number    = {3},
  pages     = {502--517},
  volume    = {26},
  doi       = {10.1109/tro.2010.2044948},
  publisher = {Institute of Electrical and Electronics Engineers (IEEE)},
}

@InProceedings{MarmoZavoli_2024_CCMfCCoLTIM,
  author    = {Marmo, Nicola and Zavoli, Alessandro},
  booktitle = {AIAA SCITECH 2024 Forum},
  title     = {Chance-constraint method for covariance control of low-thrust interplanetary missions},
  year      = {2024},
  month     = jan,
  pages     = {1--18},
  publisher = {American Institute of Aeronautics and Astronautics},
  doi       = {10.2514/6.2024-0630},
}

@Article{BlackmoreEtAl_2011_CCOPPwO,
  author    = {Blackmore, Lars and Ono, Masahiro and Williams, Brian C.},
  journal   = {IEEE Transactions on Robotics},
  title     = {Chance-Constrained Optimal Path Planning With Obstacles},
  year      = {2011},
  issn      = {1941-0468},
  month     = dec,
  number    = {6},
  pages     = {1080--1094},
  volume    = {27},
  doi       = {10.1109/tro.2011.2161160},
  publisher = {Institute of Electrical and Electronics Engineers (IEEE)},
}

@Article{NakkaChung_2023_TOoCCNSSfMPuU,
  author    = {Nakka, Yashwanth Kumar and Chung, Soon-Jo},
  journal   = {IEEE Transactions on Robotics},
  title     = {Trajectory Optimization of Chance-Constrained Nonlinear Stochastic Systems for Motion Planning Under Uncertainty},
  year      = {2023},
  issn      = {1941-0468},
  month     = feb,
  number    = {1},
  pages     = {203--222},
  volume    = {39},
  doi       = {10.1109/tro.2022.3197072},
  publisher = {Institute of Electrical and Electronics Engineers (IEEE)},
}

@Article{BoutonnetEtAl_2024_DtJT,
  author    = {Boutonnet, A. and Langevin, Y. and Erd, C.},
  journal   = {Space Science Reviews},
  title     = {Designing the JUICE Trajectory},
  year      = {2024},
  issn      = {1572-9672},
  month     = sep,
  number    = {6},
  volume    = {220},
  doi       = {10.1007/s11214-024-01093-y},
  publisher = {Springer Science and Business Media LLC},
}

@Book{AbramowitzStegun_1964_HoMFwFGaMT,
  author      = {Abramowitz, Milton and Stegun, Irene A.},
  publisher   = {Dover},
  title       = {Handbook of Mathematical Functions with Formulas, Graphs, and Mathematical Tables},
  year        = {1964},
  address     = {New York},
  edition     = {ninth Dover printing, tenth GPO printing},
  description = {BibTeX - Wikipedia, the free encyclopedia},
}

@Article{TopputoBelbruno_2015_ETwBC,
  author    = {Topputo, F. and Belbruno, E.},
  journal   = {Celestial Mechanics and Dynamical Astronomy},
  title     = {Earth–Mars transfers with ballistic capture},
  year      = {2015},
  issn      = {1572-9478},
  month     = feb,
  number    = {4},
  pages     = {329--346},
  volume    = {121},
  doi       = {10.1007/s10569-015-9605-8},
  publisher = {Springer Science and Business Media LLC},
}

@Article{Zakai_1969_OtOFoDP,
  author    = {Zakai, Moshe},
  journal   = {Zeitschrift für Wahrscheinlichkeitstheorie und Verwandte Gebiete},
  title     = {On the optimal filtering of diffusion processes},
  year      = {1969},
  issn      = {1432-2064},
  number    = {3},
  pages     = {230--243},
  volume    = {11},
  doi       = {10.1007/bf00536382},
  publisher = {Springer Science and Business Media LLC},
}

@InBook{YongZhou_1999_DPaHE,
  author    = {Yong, Jiongmin and Zhou, Xun Yu},
  chapter   = {4},
  pages     = {157--215},
  publisher = {Springer New York},
  title     = {Dynamic Programming and HJB Equations},
  year      = {1999},
  isbn      = {9781461214663},
  booktitle = {Stochastic Controls},
  doi       = {10.1007/978-1-4612-1466-3_4},
}

@InProceedings{JulierUhlmann_1997_NEotKFtNS,
  author    = {Julier, Simon J. and Uhlmann, Jeffrey K.},
  booktitle = {Signal Processing, Sensor Fusion, and Target Recognition VI},
  title     = {New extension of the Kalman filter to nonlinear systems},
  year      = {1997},
  editor    = {Kadar, Ivan},
  month     = jul,
  pages     = {1--12},
  publisher = {SPIE},
  doi       = {10.1117/12.280797},
  issn      = {0277-786X},
}

@Article{GrecoEtAl_2022_RSTDUBOC,
  author    = {Greco, Cristian and Campagnola, Stefano and Vasile, Massimiliano},
  journal   = {Journal of Guidance, Control, and Dynamics},
  title     = {Robust Space Trajectory Design Using Belief Optimal Control},
  year      = {2022},
  issn      = {1533-3884},
  month     = jun,
  number    = {6},
  pages     = {1060--1077},
  volume    = {45},
  doi       = {10.2514/1.g005704},
  publisher = {American Institute of Aeronautics and Astronautics (AIAA)},
}

@InProceedings{MarmoEtAl_2023_AHMSAfCCoIMwNE,
  author    = {Marmo, Nicola and Zavoli, Alessandro and Ozaki, Naoya and Kawakatsu, Yasuhiro},
  booktitle = {2023 AAS/AIAA Space Flight Mechanics Meeting},
  title     = {A hybrid multiple-shooting approach for covariance control of interplanetary missions with navigation errors},
  year      = {2023},
  month     = {01},
  number    = {AAS 23-208},
  pages     = {1--20},
}

@Article{OguriMcMahon_2022_SPVfRLTTDuU,
  author    = {Oguri, Kenshiro and McMahon, Jay W.},
  journal   = {Journal of Guidance, Control, and Dynamics},
  title     = {Stochastic Primer Vector for Robust Low-Thrust Trajectory Design Under Uncertainty},
  year      = {2022},
  issn      = {1533-3884},
  month     = jan,
  number    = {1},
  pages     = {84--102},
  volume    = {45},
  doi       = {10.2514/1.g005970},
  publisher = {American Institute of Aeronautics and Astronautics (AIAA)},
}

@Unpublished{CalebEtAl_2025_APBCSfFOLTTO,
  author    = {Caleb, Thomas and Armellin, Roberto and Boone, Spencer and Lizy-Destrez, Stéphanie},
  note      = {Preprint},
  title     = {Taylor polynomial-based constrained solver for fuel-optimal low-thrust trajectory optimisation},
  year      = {2025},
  copyright = {Creative Commons Attribution Non Commercial No Derivatives 4.0 International},
  doi       = {10.48550/ARXIV.2502.00398},
  publisher = {arXiv},
  url       = {https://doi.org/10.48550/ARXIV.2502.00398},
}

@Unpublished{CalebEtAl_2025_CCTaFREfSTO,
  author    = {Caleb, Thomas and Armellin, Roberto and Lizy-Destrez, Stéphanie},
  title     = {Chance constraints transcription and failure risk estimation for stochastic trajectory optimisation},
  year      = {2025},
  copyright = {Creative Commons Attribution Non Commercial No Derivatives 4.0 International},
  doi       = {10.48550/ARXIV.2502.15949},
  publisher = {arXiv},
  url       = {https://doi.org/10.48550/arXiv.2502.15949},
}

@Article{NgangaWensing_2021_ASODDPfRBS,
  author    = {Nganga, John N. and Wensing, Patrick M.},
  journal   = {IEEE Robotics and Automation Letters},
  title     = {Accelerating Second-Order Differential Dynamic Programming for Rigid-Body Systems},
  year      = {2021},
  issn      = {2377-3774},
  month     = oct,
  number    = {4},
  pages     = {7659--7666},
  volume    = {6},
  doi       = {10.1109/lra.2021.3098928},
  publisher = {Institute of Electrical and Electronics Engineers (IEEE)},
}

@InProceedings{OguriLantoine_2022_SSCPfRLTTDuU,
  author    = {Oguri, Kenshiro and Lantoine, Gregory},
  booktitle = {2022 AAS/AIAA Astrodynamics Specialist Conference},
  title     = {Stochastic Sequential Convex Programming for Robust Low-thrust Trajectory Design under Uncertainty},
  year      = {2022},
  month     = {08},
  number    = {AAS 22-708},
  pages     = {1--20},
}
